\documentclass[12pt]{amsart}
\usepackage{amsmath,amscd,amssymb,amsfonts}
\usepackage[all]{xy}

\newcommand{\ver}{October 19, 2006; v6. November 21, 2008}
\newcommand{\bC}{{\mathbb C}}
\newcommand{\bD}{{\mathbb D}}

\newcommand{\bH}{{\mathbb H}}

\newcommand{\bN}{{\mathbb N}}
\newcommand{\bP}{{\mathbb P}}
\newcommand{\bR}{{\mathbb R}}
\newcommand{\bQ}{{\mathbb Q}}
\newcommand{\bZ}{{\mathbb Z}}
\newcommand{\cA}{{\mathcal A}}
\newcommand{\cB}{{\mathcal B}}

\newcommand{\cG}{{\mathcal G}}

\newcommand{\cL}{{\mathcal L}}
\newcommand{\cM}{{\mathcal M}}
\newcommand{\cJ}{{\mathcal J}}
\newcommand{\cO}{{\mathcal O}}
\newcommand{\cP}{{\mathcal P}}
\newcommand{\cT}{{\mathcal T}}
\newcommand{\cV}{{\mathcal V}}
\newcommand{\cW}{{\mathcal W}}
\newcommand{\cU}{{\mathcal U}}

\newcommand{\rank}{{\rm rank}}
\newcommand{\Hom}{{\rm{Hom}}}
\newcommand{\Pic}{{\rm{Pic}}}

\newcommand{\wti}{\widetilde}
\newcommand{\tX}{\widetilde{X}}
\newcommand{\tZ}{\widetilde{Z}}

\newcommand{\Gr}{\text{\rm Gr}}

\newcommand{\ra}{\rightarrow}
\newcommand{\rndown}[1] {\llcorner {#1} \lrcorner}

\theoremstyle{plain}
\newtheorem{thm}{Theorem}[section]
\newtheorem{cor}[thm]{Corollary}

\newtheorem{lem}[thm]{Lemma}
\newtheorem{prop}[thm]{Proposition}
\newtheorem{que}[thm]{Question}

\theoremstyle{definition}
\newtheorem{df}[thm]{Definition}
\newtheorem{rem}[thm]{Remark}
\newtheorem{example}[thm]{Example}

\title[Local systems and multiplier ideals]{Unitary local systems,
multiplier ideals, and polynomial periodicity of Hodge numbers}
\author{Nero Budur}
\address{Department of Mathematics, The
University of Notre Dame, IN 46556, USA} \email{nbudur@nd.edu}

\date{\ver}

\thanks {The author was supported in part by the NSF grant DMS-0700360 and by the Institute of Advanced Study.}

\begin{document}

\maketitle

\begin{abstract}
The space of unitary local systems of rank one on the complement
of an arbitrary divisor in a complex projective algebraic variety
can be described in terms of parabolic line bundles.  We show that
multiplier ideals provide  natural stratifications of this space.
We prove a structure theorem for these stratifications in terms of
complex tori and convex rational polytopes, generalizing to the
quasi-projective case results of Green-Lazarsfeld and Simpson. As an application we show the polynomial periodicity of Hodge numbers $h^{q,0}$ of congruence
covers in any dimension, generalizing results of E. Hironaka and
Sakuma. We extend the structure theorem and polynomial periodicity
to the setting of cohomology of unitary local systems. In
particular, we obtain a generalization of the polynomial
periodicity of Betti numbers of unbranched congruence covers due
to Sarnak-Adams. We derive a geometric characterization of finite
abelian covers, which recovers the classic one and the one of
Pardini. We use this, for example, to prove a conjecture of
Libgober about Hodge numbers of abelian covers.
\end{abstract}


\section{Introduction}

\noindent{\bf Motivation.} This work arose out of an effort to
understand global invariants of singularities of divisors in
nonsingular complex projective varieties. One of the most
versatile tools that measures the complexity of singularities is
the multiplier ideal, developed by Nadel \cite{Nad},
Esnault-Viehweg \cite{EV2}, Demailly \cite{Dem}, Ein-Lazarsfeld
\cite{EL}, Siu \cite{AS} (see \cite{La} for more on them). Their
nature and power still remains  somewhat mysterious. Locally, there
is a better understanding of multiplier ideals through different
interpretations (resolution of singularities, plurisubharmonic
functions, positive characteristic methods, D-modules)  and
through connections with many other different notions
(log-canonical threshold, jet-schemes, monodromy on Milnor fibers, action of the Galois group on the homology of infinite abelian covers of complements to isolated non-normal crossings divisors,
Bernstein-Sato polynomial, Hodge filtration of some Hodge modules,
etc.). However, some of the most striking applications of
multiplier ideals are to global problems: Minimal Model Program,
invariance of plurigenera, linear series, etc. (see \cite{La}).
The key feature that makes the multiplier ideals so powerful is
the vanishing of spaces of the type
\begin{align}\label{moti}
H^q(X,\omega_X\otimes L\otimes \cJ(D))
\end{align}
where $X$ is a nonsingular complex projective variety, $\omega_X$
the canonical sheaf, $\cJ(D)$ the multiplier ideal of an effective
$\bQ$-divisor $D$ on $X$, and $L$ is a line bundle such that $L-D$
enjoys some positivity. In this paper we adopt an opposite point
of view. Namely, when the support of $D$ is fixed and $L-D$ is not
positive anymore we regard the dimensions of the spaces
(\ref{moti}) as global invariants of the singularities of the
support of $D$. This paper is focused on these invariants and
tries to develop a natural setting for them. Not surprisingly, in
light of work of Esnault-Viehweg \cite{EV1}, \cite{EV2}, the
natural setting we found most natural is the theory of local
systems of rank one. As a by-product of work aimed at multiplier
ideals, we are able to prove a few other general facts about
local systems.

\medskip

\noindent {\bf Local systems of rank one.} Let $X$ be a
nonsingular projective variety over $\bC$. Let $M$ be the moduli
space of complex local systems of rank one on $X$. Then $M$ has
three interpretations, that is it has three different complex
structures on the same real analytic group (see \cite{Si}) :
\begin{align*}
M_B &=\Hom (\pi _1(X),\bC^*)=\Hom (H_1(X,\bZ),\bC
^*)=H^1(X,\bC^*),\\
M_{DR} &=\{ (L,\nabla) : L\in \Pic(X),\ \nabla:L\ra
L\otimes\Omega_X^1\ \rm{integrable\ connection}\},\\
M_{Dol}&=\{(E,\theta) : \rm{Higgs\ line\
bundle}\}=\Pic^\tau(X)\times H^0(X,\Omega_X^1),
\end{align*}
where $\Pic^\tau(X)=\ker [c_1:\Pic(X)\ra H^2(X,\bR)]$. Let
$U\subset M$ be the unitary local systems. We have $U_B
=\Hom(H_1(X,\bZ), S^1)$, where $S^1\subset \bC^*$ is the unit
circle, and $U_{Dol}=\Pic^\tau(X)$. $U_B$ is a totally real
subgroup of $M_B$, whereas $U_{Dol}$ is a finite disjoint union of
copies of the Picard variety $\Pic^0(X)$. For example, for a curve
of genus $g$, we have $U_B=(S^1)^{2g}$ and $U_{Dol}={\rm
{Jac}}(X)$, the Jacobian of $X$.

\smallskip

\noindent{\it Notation.}  Let $D=(D_i)_{i\in S}$ be a finite tuple
of distinct irreducible and reduced codimension one subvarieties
of $X$. We use the following notation:
\begin{align*}
\alpha+\beta&=(\alpha_i+\beta_i)_{i\in S} & \alpha\cdot
D&=\sum_{i\in S}\alpha_i\cdot D_i\\
\rndown{\alpha}&=(\rndown{\alpha_i})_{i\in S} & \alpha\cdot
[D]&=\sum_{i\in S}\alpha_i\cdot [D_i]\\
\{\alpha\}&=(\{\alpha_i\})_{i\in S} & \mu^*(\alpha\cdot D) &=\sum_{i\in S}\alpha_i\cdot \mu^*D_i
\end{align*}
where $\alpha,\beta\in \bR^S$, $[D_i]$ are the  cohomology classes
of $D_i$ (it will be clear from context which cohomology),
$\mu:Z\ra X$ is a map of varieties, $\rndown{.}$ is the round-down
function, $\{.\}$ is the fractional part. We also denote by $D$
the divisor $\cup_{i\in S}D_i$. By $H_*(X)$ ( $H^*(X)$) we always
mean the (co)homology  with integral coefficients. We will use the
term {\it polytope} to mean a usual polytope that does not
necessarily include all its faces.

\begin{df} (a) Define the set of {\it realizations of boundaries of $X$ on $D$} as
$$
\Pic^\tau(X,D):=\left\{ (L,\alpha) \in\Pic(X)\times [0,1)^S :
c_1(L)=\alpha\cdot [D]\ \in H^2(X,\bR)\right\}.
$$
(b) Define the set of {\it realizable boundaries of $X$ on $D$} as
$$
B(X,D):=\left\{ \alpha\in [0,1)^S : c_1(L)=\alpha\cdot [D]\in
H^2(X,\bR) \text{ for some }L\in\Pic(X)\right\}.
$$
\end{df}

It is  common in the literature to call the realizations of
boundaries on $X$ {\it parabolic line bundles} of parabolic degree
0. They were introduced by Seshadri in \cite{Se}. The
terminology is influenced by the use of the term boundary in
birational geometry. However, from our point of view it is natural
to exclude the boundaries for which some coefficient is 1.

\medskip

The set of realizations of boundaries $\Pic^\tau(X,D)$ is an
abelian group under the operation
$$(L_,\alpha)\cdot
(L',\alpha')=\left (L\otimes L'\otimes\cO_X(-
\rndown{\alpha+\alpha'}\cdot D), \{\alpha+\alpha'\}\right ).$$ The
operation of taking the fractional part of the sum induces a group
structure on the set of realizable boundaries $B(X,D)$.

\begin{thm}\label{ls} Let $X$ be a nonsingular
complex projective variety, $D$ a divisor on $X$, and let $U=X-D$.
There is a natural canonical group isomorphism
$$RH: \Pic^\tau(X,D)\xrightarrow{\sim}\Hom (H_1(U),S^1)$$
between realizations
of boundaries of $X$ on $D$ and unitary local systems of rank one on $U$.
\end{thm}
The original result is the case of $X$ being a projective curve
and for higher rank local systems and is due to Mehta-Seshadri
\cite{MS} (announced in \cite{Se}). In this paper, for simplicity
we only consider rank one local systems. The theorem is not new,
it is known to differential geometers at least in the simple
normal crossings case, but we could not find a reference for a
simple self-contained proof. The isomorphism $RH$ follows for
example from T. Mochizuki's  extension in \cite{MoII} of the
Kobayashi-Hitchin correspondence from \cite{Si1} to the open case
and probably from other partial results leading to this extension
(see Jost-Zuo \cite{JS}, Simpson \cite{Si2}, Biquard \cite{Bi}, J.
Li \cite{LiJ}, Steer-Wren \cite{SW}). We will reproduce a
self-contained proof of Theorem \ref{ls}, showed to us by T.
Mochizuki, which is fairly elementary, requiring only standard
tools from differential geometry. This shortens the exposition
from an earlier preprint in which we derived Theorem \ref{ls} from
the powerful Theorem 1.1 of \cite{MoII}. We show the isomorphism
$RH$ is natural (Proposition \ref{comp}). For $N> 1$, the
character group of $H_1(U,\bZ_N)$ is the $N$-torsion part of
$\Pic^\tau(X,D)$ (Lemma \ref{torsion}). The group of realizations
of boundaries depends only on the complement $U$; see Proposition
\ref{comparison} for the explicit isomorphism between the groups
of realizations of boundaries on a pair $(X,D)$ and on a
resolution $(Z,E)$ isomorphic above $U$. On $\Pic^\tau (X,D)$
there is a mixture of real and complex structures: $\Pic^\tau
(X,D)$ maps onto $B(X,D)\subset\bR^S$ with fibers $\Pic^\tau(X)$
consisting of finitely many disjoint copies of the complex torus
$\Pic^0(X)$. The purely real part $B(X,D)$ splits canonically into
a finite disjoint union of convex rational polytopes in $\bR^S$
(see Section \ref{ns}).

\medskip

\noindent {\bf Canonical stratifications.} The space of unitary
local systems of rank one on a quasi-projective nonsingular
variety $U$ has canonical stratifications. Given $U$, $X$, and $D$
as above define
$$V^q_i(X,D):=\{ (L,\alpha)\in \Pic^\tau(X,D) : h^q(X, \omega_X\otimes
L\otimes\cJ (\alpha\cdot D))\ge i \},$$ where $\cJ (\alpha\cdot
D)$ denotes the multiplier ideal of the $\bR$-divisor $\alpha\cdot
D$. Under the isomorphism of Theorem \ref{ls}, $V^q_i(X,D)$
defines a subset $V^q_i(U)$ of the space of unitary local systems
of rank one on $U$ independent of $X$ and $D$ (Lemma \ref{indep}).
$V^q_i(U)$ can also be interpreted via the initial piece of the
Hodge filtration (see \cite{Ti}) on the cohomology of $U$ with
unitary local systems as coefficients (Proposition \ref{timm}).

\medskip

We prove a structure result about these canonical sets. The
original theorem when $D=0$ is due to Green-Lazarsfeld \cite{GL1},
\cite{GL2}. For different proofs see also Arapura \cite{Ar1},
Pink-Roessler \cite{PR}, Simpson \cite{Si}. The most general
theorem over projective varieties is due to Simpson \cite{Si}.
Over quasi-projective varieties, Arapura \cite{Ar2} has a similar
result for the canonical sets defined via the cohomology of the
local systems under the assumptions of $D$ having normal crossings
and $H^1(X,\bC)=0$.  Using Simpson's version and Theorem \ref{ls},
we give a generalization in a different direction  to
quasi-projective varieties which holds without additional
hypotheses.

\begin{thm}\label{gv} With the notation of Theorem \ref{ls}, there
exists a decomposition of the set of realizable boundaries
$B(X,D)\subset\bR^S$ into
 a finite number of  rational convex polytopes $P$  such that for every $q$ and $i$ the subset
$V^q_i(X,D)$ of $\Pic^\tau(X,D)$ is a finite union of sets of the
form $P\times \cU$, where $\cU$ is a torsion translate of a
complex subtorus of $\Pic^\tau(X)$. Any intersection of sets
$P\times \cU$ is also of this form. Pointwise, the subset of
$V^q_i(X,D)$ corresponding to $P\times \cU$ consists of the
realizations $(L+M,\alpha)$ with $\alpha\in P$ and $M\in \cU$, for
some line bundle $L$ depending on $P$ which can be chosen such
that $(L,\alpha)$ is torsion for some $\alpha\in P$.
\end{thm}
The decomposition into polytopes of $B(X,D)$ is not the canonical
one mentioned before, but it is any of the ones induced by
log-resolutions of $(X,D)$. When $D$ is a simple normal crossing
divisor, the decomposition of $B(X,D)$ is the canonical one.

\medskip

The proof of Theorem \ref{gv} extends to cover a more general
result regarding loci given by dimensions of pieces of the Hodge
filtration on cohomology of unitary local systems on $U$. Define
$$W^{p,q}_i(U):=\{\ \cV \in \Hom(H_1(U),S^1)\ |\ \dim
\Gr_F^pH^{p+q}(U,\cV)\ge i\ \}.$$ Then $V^q_i(U)=W^{0,n-q}_i(U)^\vee:=\{\cV^\vee\ |\ \cV\in W_i^{0,n-q}(U)\}$, see Proposition \ref{timm}.
More generally, let
$$W^{p,q}_i(U,\cW)=\{\ \cV \in \Hom(H_1(U),S^1)\ |\ \dim
\Gr_F^pH^{p+q}(U,\cW\otimes\cV)\ge i\ \}$$ for a fixed unitary
local system $\cW$ on $U$ of arbitrary rank. In this notation,
$W^{p,q}_i(U)=W^{p,q}_i(U,\bC_U)$.

\begin{thm}\label{gvgen} The statement of Theorem \ref{gv} holds
for the image of $W^{p,q}_i(U,\cW)$ in $\Pic^\tau(X,D)$ under the
assumption that either $\cW$ is the restriction of a local system
from a nonsingular compactification of $U$, or there are no
nontrivial unitary rank one local systems on a such
compactification.
\end{thm}
The decomposition of $B(X,D)$ into polytopes needed for Theorem
\ref{gvgen} is in general a refinement of the decomposition needed
for Theorem \ref{gv} and is given by Remark \ref{canonical
extension rank one}-(d).

\begin{cor}\label{corollary_to_general case} With the assumptions of Theorem \ref{gvgen}, the statement of
Theorem \ref{gv} holds for the image of
$$\{ \cV\in\Hom(H_1(U),S^1)\ |\ \dim H^m(U,\cW\otimes\cV)\ge i\ \}$$
in $\Pic^\tau(X,D)$.
\end{cor}

When $\cW=\bC_U$ and $H^1(X,\bC)=0$, the loci of Corollary
\ref{corollary_to_general case} (without the restriction to the unitary case) are considered in the
first Theorem of \cite{Ar2} in terms of maps from $U$ to products complex tori. Our method does not allow us to replace in Corollary \ref{corollary_to_general case} $S^1$ by
$\bC^*$ (i.e. unitary by arbitrary). One could ask if the polytopes in the decomposition from Corollary \ref{corollary_to_general case}, under the embedding $B(X,D)\subset (S^1)^S$, are translates of real subtori $(S^1)^r$. Our method does not answer this question.

\medskip

One problem with removing the assumption on
$\cW$ in Theorem \ref{gvgen} and Corollary
\ref{corollary_to_general case} is that we do not know the answer
of the following question, for example (see Remark \ref{on_IH}):

\begin{que}\label{qgen} With the notation as above, fix integers $m$ and $i$.
Suppose $\cW$ is a unitary local system on $U$. Is the image of
the set
$$\{\ \cV\in\Hom(H_1(X),S^1)\ |\ \dim
IH^m(X,\cW\otimes \cV_{|U})\ge i\ \}$$ in $\Pic^\tau(X)$ a finite
union of torsion translates of complex subtori of $\Pic^\tau(X)$ ?
\end{que}

{\it Connection with other recent works.} (a) Professor Libgober has
kindly informed us that for the case when $q\ge 1$, $D$ having
isolated non-normal crossings singularities, the irreducible
components $D_i$ being ample, and $H_1(X)=0$, the polytopes in the
description of $V^q_i$ are determined (see \cite{Li2}) by the
local polytopes of quasiadjunction introduced in \cite{Lib4},
extending the case $X=\bP^2$ of \cite{Lib}. In \cite{Li2}, under
the same assumptions, the sets $V^q_i$ were related to the
homotopy groups $\pi_q(U)$. See \cite{Lib}, \cite{Li2}, \cite{Li}
for properties and applications of local and global polytopes of
quasiadjunction. Even if the end result is the same as ours in
this special case, the methods are different (see Example
\ref{ex1ct1}). We do not know how to extend to method of polytopes
of quasiadjunction to recover the more general case we consider.
For a streamlined version of the method of this article for a
special case of Theorem \ref{gvgen} see Theorem 2.1 of the recent
preprint \cite{Li5} which appeared before this revision was made.

(b) The size of sets similar-looking to $V^q_i$ have been considered
from the point of view of the Fourier-Mukai transform by
Pareschi-Popa \cite{PP}. For example, for $X$ a smooth projective variety of dimension $n$ and Albanese dimension $n-k$, \cite{PP}-Theorem A gives that the codimension in $\Pic^\tau(X)$ of $V^q_1(L,\alpha)$ (see Lemma \ref{si2}) is $\ge q-k+\kappa$, where $\kappa$ is the Iitaka dimension of $L-\alpha\cdot D$ along a generic fiber of the Albanese map of $X$. If, in the Theorem \ref{gv} for the decomposition of $V^q_1(X,D)$, a set $P\times\cU$ consists of realizations $\{(L+M,\alpha)\}$, then the codimension of $\cU$ equals the codimension of $V^q_1(L,\alpha)$ for some $\alpha\in P$ (cf. proof of Theorem \ref{gv}).

\medskip

\noindent {\bf Congruence covers and polynomial periodicity of
Hodge numbers.} As an application of the above results, we use
counting of lattice points in rational convex polytopes to show
polynomial periodicity of the Hodge numbers $h^{q,0}$ of
congruence covers and, more generally, of the Hodge ranks
$h^{p,q}$ (see definition below) of the cohomology of the
pullbacks of a unitary local system to the unbranched congruence
covers.

\begin{df} A function $f$ on the set $\{1,2,3,\ldots\}$ is a {\it quasi-polynomial} if there exists a natural number $M$ and polynomials $f_i(x)\in\bQ[x]$ for $1\le i\le M$ such that $f(N)=f_i(N)$ if $N\equiv i\; ( \text{mod}\; M)$.
\end{df}

Let $X$ be a
nonsingular complex projective variety, $D=\cup_{i\in S}D_i$ a
divisor on $X$ with irreducible components $D_i$, and let $U=X-D$.
The canonical surjections
$$H_1(U,\bZ)\ra H_1(U,\bZ_N),$$
for $N>1$, give canonical normal abelian covers $X_N\ra X$ called
congruence covers. Let $h^q(N)$ denote the Hodge numbers
$h^{q,0}=h^{0,q}$ of any nonsingular model of $X_N$. These are
birational invariants hence $h^q(N)$ is well-defined.

\begin{thm}\label{pp} With the notation as above, for every $q$ the
function $h^q(N)$ is a quasi-polynomial in $N$.
\end{thm}

The case when $X$ is a surface was done by E. Hironaka \cite{Hi1},
\cite{Hi2}, and Sakuma \cite{Sa}, generalizing previous results of
Zariski, Libgober, Vacquie (see \cite{Hi1} for the history). Even the case of $X=\bP^n$ for $n>2$ is new to our knowledge.
The original result is of Sarnak-Adams \cite{SA} who proved in
every dimension the polynomial periodicity of the  Betti
numbers of the unbranched congruence covers $U_N\ra U$ (see also \cite{CO}). For a survey of applications and related results in the case of hyperplane arrangements see \cite{Su}.

\medskip

The proof of Theorem \ref{pp}, adapted to the situation of Theorem
\ref{gvgen}, yields the following generalization. Let $g_N:U_N\ra
U$ be the unbranched congruence cover and fix a unitary local
system $\cW$ on $U$. Define the {\it Hodge rank} $h^{p,q}_\cW(N)$
as the dimension of $\Gr_F^pH^{p+q}(U_N,g_N^*\cW)$. This
definition should not be confused with the usual definition of
Hodge numbers of a mixed Hodge structure which involves the weight
filtration as well. Typically $H^{k}(U_N,g_N^*\cW)$ does not have
weight $k$ and it is not pure. However, we do not worry about the
weight filtration in this article. If $\cW=\bC_U$ and $X_N$ is a
nonsingular compactification of $U_N$ with a complement $E_N$
which is a simple normal crossings divisor, then
$h^{p,q}_\cW(N)=h^q(X_N,\Omega_{X_N}^p (\log E_N))$ and
$h^{0,q}_\cW(N)=h^q(N)$ as in Theorem \ref{pp}.

\medskip

We assume that $\cW$ is the restriction of a local system from a
nonsingular compactification of $U$, or that there are no
nontrivial unitary rank one local systems on a such
compactification. The following generalizes the result of
Sarnak-Adams \cite{SA} on the Betti numbers of $U_N$:

\begin{thm}\label{ppgen} For every pair of integers $p$ and $q$, the Hodge rank function
$h^{p,q}_\cW(N)$ is a quasi-polynomial in $N$.
\end{thm}

\medskip

\noindent {\bf Finite abelian covers.}  The subject of finite
abelian covers is well-known and has been studied by many people.
To our knowledge, the geometric characterization of abelian covers
has two versions: a classical one (\cite{Za}, \cite{Na}) and the
one of \cite{Pa}. The isomorphism of Theorem \ref{ls} recovers
both characterizations offering a common generalized approach.
This does not seem to be well-known to algebraic geometers since
as a consequence we are able to prove some results about abelian
covers (Corollaries \ref{ss2} and \ref{ss3}) which, to our
knowledge, appear in their published version only under some
additional hypotheses. We employ these results in deriving the
theorems mentioned earlier. From Theorem \ref{ls} we derive the
following geometric characterization, essentially proven in
\cite{Na}.

\begin{cor}\label{fc} Let $X$ be a nonsingular complex projective
variety, $D=\cup_{i\in S}D_i$ a divisor on $X$ with irreducible
components $D_i$, and let $U=X-D$. Let $G$ be a finite abelian group. The  equivalence classes of
normal $G$-covers (Definition \ref{gcovering}) of $X$ unramified  above $U$ are into
one-to-one correspondence with the subgroups
$G^*\subset\Pic^\tau(X,D)$.
\end{cor}

Here $G^*=\Hom(G,\bC^*)$ is the dual group of $G$. We show that
Corollary \ref{fc} also recovers via the isomorphism $RH$ the geometric characterization of
finite abelian covers due to Pardini \cite{Pa}. In
particular, as expected, one has the following:

\begin{cor}\label{ss} With the notation as in Corollary \ref{fc}, let
$\pi:\tX\ra X$ be a normal $G$-cover of $X$ unramified above
$U$ corresponding to an inclusion $G^*=\{(L_\chi ,\alpha_\chi ):
\chi\in G^*\}\subset\Pic^\tau(X,D)$. Then
$$\pi_*\cO_{\tX}=\bigoplus_{\chi\in G^*}L_\chi^{-1},$$
where $G$ acts on $L_\chi^{-1}$ via the character $\chi$.
\end{cor}

The following computation was done by Esnault-Viehweg \cite{EV2} in
the case when $G$ is cyclic.

\begin{cor}\label{ss2} With the notation as in Corollary \ref{ss}, suppose
$\mu:Z\ra X$ is a log resolution of $(X,D)$ which is an
isomorphism above $U$. Let $\rho:\tZ\ra Z$ be the corresponding
normal $G$-cover unramified above $U$. Then we have the
eigensheaf decomposition
$$\rho_*\cO_{\tZ}=\bigoplus_{\chi\in
G^*}\mu^*L_\chi^{-1}\otimes\cO_Z(\rndown{\mu^*(\alpha_{\chi}\cdot
D)}),$$ where the round-down of a divisor means rounding-down
of the coefficients of its irreducible components.
\end{cor}

The following computation was done by Libgober (see \cite{Li2}) in the case $X=\bP^N$ and proposed as Conjecture 6.4 in \cite{Li2} for the case of $H_1(X)=0$, $D$ having ample irreducible components and isolated non-normal crossings, and $q=n$:

\begin{cor}\label{ss3} With the notation as in Corollary \ref{ss}, let
$H^0(Y, \Omega_{Y}^q)$ denote the space of global  q-forms on a
nonsingular model $Y$ of $\tX$. Then
$$H^0(Y, \Omega_{Y}^q)\cong\bigoplus_{\chi\in
G^*}H^{n-q}(X,\omega_X\otimes L_\chi\otimes\cJ(\alpha_{\chi}\cdot
D)),$$ where $n$ is the dimension of $X$.
\end{cor}

\medskip

\noindent{\bf Layout.} In section \ref{misection} we review
multiplier ideals,  local systems, and the details of analysis
needed for the next section.  In section \ref{ulsrb} we give a
self-contained proof of Theorem \ref{ls}. In section \ref{ns} we
show the naturality of this isomorphism and describe the structure
of the spaces of realizations of boundaries and of realizable
boundaries. In section \ref{facs} we derive the geometric
characterization and computations for finite covers described in
the Introduction. In section \ref{canstrat} we prove Theorem
\ref{gv}. In section \ref{congcov} we prove Theorem \ref{pp}. In
the last section we prove generalizations in terms of local
systems of the above results, namely Theorem \ref{gvgen}, Corollary \ref{corollary_to_general case}, and
Theorem \ref{ppgen}.

\medskip

\noindent{\bf Acknowledgement.} We thank  L. Ein, L. Maxim, J.
Song, the referees, and especially M. Saito for helpful comments and
discussions. We thank T. Mochizuki for telling us how to give a
self-contained proof of Theorem \ref{ls}, shortening the
exposition of an earlier preprint. We thank the Institute of
Advanced Study in Princeton for their hospitality.

\section{Basic notions}\label{misection}

In this paper, all algebraic varieties are defined over the field of complex numbers. The underlying analytic variety of an algebraic variety $X$ will be denoted by the same notation. This section recalls the definition of the basic objects we work with in this paper and it can be skipped by the knowledgeable reader.

\medskip
\noindent{\bf Multiplier ideals.} We review the definition of multiplier ideals and the local vanishing theorem (\cite{La}). For a nonsingular variety $X$, let $\omega_X$ denote the canonical sheaf.  Let $\mu:Y\ra X$ be a proper birational morphism. The exceptional set of $\mu$ is denoted by $Ex(\mu)$. This is the set of points $\{y\in Y\}$ where $\mu$ is not biregular. We say that $\mu$ is a resolution if $Y$ is nonsingular. Let $D$ be an effective $\bR$-divisor on $X$. We say that a resolution $\mu$ is a {\it log resolution} of $(X,D)$ if  $\mu^{-1}(D)\cup Ex(\mu)$ is a divisor with simple normal crossings. Such a resolution always exists, by Hironaka. Denote by $\omega_{Y/X}$ the relative canonical sheaf $\omega_Y\otimes(\mu^*\omega_X)^{\otimes -1}$. If $D=\sum \alpha_iD_i$ is an
    effective $\bR$-divisor on $X$,
  where $D_i$ are the irreducible components of $D$ and $\alpha_i\in\bR$, the {round down}
  of $D$ is the integral divisor $\rndown{D}=\sum\rndown{\alpha_i}D_i$.

\begin{df} Let $D$ be an effective $\bR$-divisor on a nonsingular complex variety $X$, and fix a log resolution $\mu:Y\ra X$ of $(X,D)$. Then the {\it multiplier ideal sheaf} of $D$ is defined to be $\cJ(D)=\mu_*(\omega_{Y/X}\otimes\cO_Y(-\rndown{\mu^*D}))$.
\end{df}

This is a sheaf of ideals in $\cO_X$ independent of the log resolution. Analytically, the multiplier ideals are defined locally to consist of holomorphic functions $f$ such that $|f|^2/\prod |g_i|^{2\alpha_i}$ is locally integrable, where $g_i$ are fixed equations of the irreducible components of $D$ and $\alpha_i\in\bR_{>0}$ are the corresponding coefficients in $D$.

\begin{thm}\label{localvanishing} (Local vanishing, \cite{La}-Theorem 9.4.1) Let $D$ be an effective $\bR$-divisor on a nonsingular complex variety $X$, and $\mu:Y\ra X$ of $(X,D)$ a log resolution of $(X,D)$. Then $R^j\mu_*(\omega_{Y/X}\otimes\cO_Y(-\rndown{\mu^*D}))=0$ for $j>0$.
\end{thm}

\medskip
\noindent{\bf Metrics, connections, and curvature.} Recall the following basic facts
from differential geometry (\cite{GH}). Let $X$ be a
complex manifold and $L$ a complex line bundle on $X$. A
trivialization of $L$ over an open subset $\Omega$ of $X$ is the
same as a frame $u$ of $L$ above $\Omega$. Fix frames $u$ for $L$.
A {\it hermitian metric} $(.\; ,.)$ on $L$ is equivalent to a
$C^\infty$ real-valued positive function $h(x)=(u(x),u(x))$. One also writes $h=\exp(-2\phi)$, with $\phi\in
C^\infty(\Omega,\bR)$. A {\it connection} on $L$ is a linear map
$$\bD:\cA^0(L)\ra \cA^1(L)$$
where $\cA^p(L)$ is the sheaf of $C^\infty$ $L$-valued $p$-forms
on $X$ such that $\bD(\sigma u)=(d\sigma) u +\sigma \bD(u)$ for
$C^\infty$ functions $\sigma$ where $d$ is the usual exterior
derivative. If $L$ is a holomorphic line bundle, there is a
natural $(0,1)$-operator $\bar{\partial}$ on $\cA^0(L)$ defined
locally by $\bar{\partial}(\sigma u)=(\bar{\partial}{\sigma})u$
for a holomorphic frame $u$. In this case there exists a unique
connection $\bD$ on $L$ such that it is compatible with the
hermitian metric $h$ and with the complex structure. That is
$d(v,w)=(\bD v,w)+(v,\bD w)$ and, if $\bD=\bD'+\bD''$ is the decomposition of
$\bD$ into the $(1,0)$ and $(0,1)$ parts, then $\bD''=\bar{\partial}$.
$\bD$ is called the {\it metric connection} of $h$. If $u$ is a
holomorphic frame for $L$ and $\bD$ is the metric connection of $h$,
then
$$\bD u=\frac{\partial h}{h}\cdot u.$$ A connection on a holomorphic line bundle $L$ is a {\it holomorphic connection} if it compatible with the complex structure. A connection $\bD$ extends to operators $\bD:\cA^p(L)\ra\cA^{p+1}(L)$
by forcing $\bD(\sigma u)=(d\sigma) u+(-1)^p\sigma\wedge \bD u$ for all  $C^\infty$ $p$-forms $\sigma$. The operator $\bD^2:\cA^0(L)\ra\cA^2(L)$ is linear in the $C^\infty$ functions on  $X$ and $\bD^2u=\Theta u$ for a closed globally defined
2-form $\Theta$, called the {\it curvature} of $\bD$. If $u$ is a
holomorphic frame for $L$ and $\bD$ is the metric connection of
$h=\exp (-2\phi)$, then
\begin{equation}\label{curvature}
\frac{i}{2\pi}\Theta=\frac{i}{\pi}\partial\bar{\partial}\phi.
\end{equation}
If in addition $X$ is compact, then the real cohomology class
defined by $i(2\pi)^{-1}\Theta$ on $X$ is the same as the first
Chern class $c_1(L)$.

\medskip
\noindent{\bf Local systems.} A {\it complex local system} $\cV$
on a complex manifold $X$ is a locally constant sheaf of finite
dimensional complex vector spaces (\cite{Di}). The rank of $\cV$
is the dimension of a fiber of $\cV$. Local systems of rank one on
$X$ are equivalent with representations $H_1(X)\ra \bC^*$. Unitary
local systems of rank one on $X$ correspond to representations
$H_1(X)\ra S^1$, where $S^1$ is the unit circle in $\bC$. Let $L$
be a complex line bundle on $X$. If $\bD^2=0$ for a connection
$\bD$ on $L$, then $\bD$ is called {\it flat}. It can be showed
then that there is exactly one holomorphic structure on $L$ such
that the $(0,1)$-part of $\bD$ is compatible with it. The functor
taking a complex bundle $L$ on $X$ with a flat connection $\bD$ to
the kernel of $\bD$ defines an equivalence between the categories
of complex vector bundles with flat connections and local systems.
Equivalently, there is an equivalence between the categories of
holomorphic vector bundles with flat holomorphic connections and
local systems. Under this equivalence, a local system is unitary
if and only if the corresponding flat connection is the metric
connection for some hermitian metric (\cite{Dem3}-V). It follows
from (\ref{curvature}) that in this case the flatness condition on
$\bD$ reads, in terms of a holomorphic frame $u$ for $L$ and
metric $h=\exp (-2\phi)$, as $\partial\bar{\partial}\phi=0$.

\medskip
\noindent {\bf The canonical Deligne extension.} Let $X$ be a
compact complex manifold and  $U$ an open subset such that the
complement is a divisor with simple normal crossings. Let $L_U$ be
a holomorphic line bundle on $U$ with a holomorphic flat
connection $\bD$ defining a unitary local system on $U$. Let $D_i$
be the irreducible components of $X-U$ and let $\alpha_i\in [0,1)$
be such that the monodromy of $(L_U,\bD)$ around a general point
of $D_i$ is multiplication by $\exp(-2\pi i\alpha_i)$. There
exists an extension $(L,\bD)$ to a holomorphic line bundle with a
flat logarithmic connection on $X$ uniquely characterized by the
condition that the eigenvalue of the residue of $(L,\bD)$ around
$D_i$ is $\alpha_i$ (\cite{De}). $(L,\bD)$ is called the {\it
canonical Deligne extension} of $(L_U,\bD)$. We recall next how
the line bundle $L$ is constructed. Explicitly, let $u_i$ be a
frame of the multi-valued flat sections of $(L_U,\bD)$ around a
general point of $D_i$. We can assume that a local equation for
$D_i$ is given by the local coordinate $z_i$. Let $v_i=\exp(\log
z_i\cdot \alpha_i)\cdot u_i$. Then $v_i$ is a single-valued
holomorphic section of $L_U$. Let $L_0$ be the line bundle on
$X-W$ obtained by declaring $v_i$ to be a local holomorphic frame,
where $W$ is the closed set of singular points of $\cup_iD_i$.
Since the codimension of $W$ is $\ge 2$, $L_0$ extends uniquely to
a line bundle $L$ on $X$.
\begin{prop}\label{chern}
 (\cite{EV1}-(B.3)) With the notation as above, the first Chern class of $L$ is given by $c_1(L)=-\sum_{i\in S}\alpha_i [D_i]$ in $H^2(X,\bR)$.
\end{prop}

\section{Unitary local systems and realizations of
boundaries}\label{ulsrb}

\medskip
\noindent {\bf Proof of Theorem \ref{ls} for the simple normal
crossings case.}  With the notation of Theorem \ref{ls}, assume in
addition that $D$ is a simple normal crossings divisor. We define
a map
$$RH:\Pic^\tau(X,D)\ra \Hom (H_1(U),S^1)$$
as follows, by an argument showed to us by T. Mochizuki. This
shortens the exposition from an earlier preprint in which we
derived the map $RH$ from the powerful correspondence of Theorem
1.1 of \cite{MoII}. The responsibility for any faults of the
exposition lies with the author.

\medskip

Let $(L,\alpha)$ be a realization of a boundary. To attach a
unitary local system to $(L,\alpha)$ we will find a hermitian
metric $h$ on the restriction of $L$ to $U$ such that the
curvature of the metric connection $\bD_h$ is zero. Then, by
section \ref{misection}, $(L_{|U},\bD_h)$ gives a unitary local
system of rank one on $U$.


\medskip

Fix an open cover $\mathcal{C}$ of $X$ and trivializations above
$\Omega\in\mathcal{C}$ of $L$. Let $h_0$ be a hermitian metric on
$L$. Let $\sigma_i$ be a nonzero holomorphic global section of the
line bundle $\cO_X(D_i)$ vanishing exactly at $D_i$. Let
$|\sigma_i|:X\ra\bR$ be the norm of $\sigma_i$ for a fixed
hermitian metric $H_i$ on $\cO_X(D_i)$. Let $H_{i,U}$ be the
restriction of $H_i$ to $U$. Remark that $\sigma_
{i,U}={\sigma_i}_{|U}$ is a frame for $\cO_X(D_i)$ above $U$.
Hence we can assume $H_{i,U}=|\sigma_{i,U}|^2$. Define $h_{1,U\cap
\Omega}= h_{0,U\cap \Omega}\prod_{i\in S}|\sigma_{i,U\cap
\Omega}|^{-2\alpha_i}$, where the appearance of $U\cap\Omega$ in
the index means restriction of the corresponding function to
$U\cap\Omega$. Then $h_{1,U\cap \Omega}$ are $\bR_{>0}$-valued
smooth functions on $U\cap \Omega$. For $\Omega, \Omega
'\in\mathcal{C}$, denote by $g_{\Omega, \Omega '}$ the transition
functions for the trivializations of $L$ above
$\Omega\cap\Omega'$. Then, on $U\cap\Omega\cap\Omega'$, we have
$h_{1,U\cap\Omega'}=|g_{\Omega, \Omega '}|^2\cdot
h_{1,U\cap\Omega}$. In other words, the $h_{1,U\cap \Omega}$ patch
up to form a hermitian metric on $L_{|U}$. We call this metric
$h_{1,U}$. Remark that, unless all the $\alpha _i$ are $0$, the
same process does not define a metric on $L$.

\medskip

By (\ref{curvature}), there is an equality of $(1,1)$-forms
globally defined on $U$
\begin{align}\label{curv}
\frac{i}{2\pi}\Theta_{h_{1,U}}=\frac{i}{2\pi}\Theta_{h_{0,U}}-\sum_{i\in
S}\alpha_i\frac{i}{2\pi}\bar{\partial}\partial\log |\sigma_{i,U}|^2,
\end{align}
where $\Theta_{*}$ is the curvature of the metric connection of
the metric $*$. The closed $(1,1)$-form
$i(2\pi)^{-1}\Theta_{h_{0,U}}$ is the restriction to $U$ of the
closed $(1,1)$-form $i(2\pi)^{-1}\Theta_{h_{0}}$ which represents
the first Chern class of $L$. The closed $(1,1)$-form
$i(2\pi)^{-1}\bar{\partial}\partial\log |\sigma_{i,U}|^2$ is the
restriction to $U$ of the closed $(1,1)$-form on $X$ given by
$i(2\pi)^{-1}\Theta_{H_i}$ which represents the cohomology class
of $D_i$. Define a global $(1,1)$-form on $X$ by
$$\frac{i}{2\pi}\Theta_{h_1}=\frac{i}{2\pi}\Theta_{h_{0}}-\sum_{i\in
S}\alpha_i\frac{i}{2\pi}\Theta_{H_i}.$$ Then $\Theta_{h_{1,U}}$ is
the restriction to $U$ of $\Theta_{h_1}$. By assumption
$c_1(L)=\alpha\cdot [D]$ in $H^2(X,\bR)$. Hence the cohomology
class of $i(2\pi)^{-1}\Theta_{h_1}$ is zero. By the
$\bar{\partial}\partial$-Poincar\'{e} lemma (\cite{GH}-p.387),
$i(2\pi)^{-1}\Theta_{h_1}=i(2\pi)^{-1}\bar{\partial}\partial\chi$
for some real-valued function $\chi$. Define a new hermitian
metric on $U$ by
$$h=h_{1,U}\exp(-\chi_{|U}).$$
Then the curvature of the metric connection $\bD_h$ of $h$ is zero.
Hence $(L_{|U},\bD_h)$ is flat connection over $U$.

\medskip

The map $RH$ is well-defined. A different choice for the metric
$h_0$ gives a new metric $h'$ which differs from $h$ by
multiplication by a constant, hence $\bD_h$ and $\bD_{h'}$ define
the same local system. It is straight forward that $RH$ is a group
homomorphism.

\medskip

The monodromy of the local system attached to $(L,\alpha)$ around
a general point of $D_i$ is multiplication by $\exp(2\pi
i\alpha_i)$. To observe this let $\Omega\subset X$ be a small
neighborhood around a general point of $D_i$. We can assume that
$\Omega=\Delta^ n$, where $n$ is the dimension of $X$ and $\Delta$
is a small disc in $\bC$, with coordinates $z_1,\ldots ,z_n$. We
can further assume that $z_1$ is a local equation of
$D_{|\Omega}=(D_i)_{|\Omega}$. Let $\Omega^*$ be $\Omega-D$, which
we can assume to be $\Delta_1^*\times\Delta^{n-1} $, where
$\Delta_1^*$ is the punctured disc in $\bC$ with coordinate $z_1$.
Then $L_{|\Omega^*}=\cO_{\Omega^*}\cdot u$, where $u$ is a fixed
holomorphic frame for $L$ above $\Omega^*$ which we can assume to
be orthonormal for $h_0\prod_{j\ne
i}H_{j}^{-\alpha_j}\exp(-\chi)$. Then $h=|z_1|^{-2\alpha_i}$. The
metric connection of $h$ is given by
$$\bD u=\frac{\partial h}{h}\cdot u=(-\alpha_i)\frac{dz_1}{z_1}\cdot
u.$$ The flat sections of $L_{|\Omega^*}$ are of the form $\sigma
u$ such that
\begin{align}\label{kercomp}
0=\bD(\sigma u)=\left
(d\sigma-\alpha_i\sigma\frac{dz_1}{z_1}\right )u.
\end{align}
Thus a multi-valued frame for the local system corresponding to
the solution space of (\ref{kercomp}) is $u_i=\exp(\log z_1 \cdot
\alpha_i)\cdot u$ (see also \cite{De}-II.1.17.1). Let $T_i$ be the
monodromy around $D_i$ of this local system. Then
$$T_iu_i=\exp((\log z_1 +2\pi i) \cdot
\alpha_i)\cdot u=\exp(2\pi i\alpha_i) u_i,$$ which is what we
wanted to show.

\medskip

The map $RH$ is injective. Indeed, if $\bC_U$ is the resulting
local system, its only extension to $X$ as a local system is
$\bC_X$. On the other hand, since the monodromy around the $D_i$'s
is trivial, we must have $\alpha_i=0$ for all $i\in S$. Thus
$h=h_0\exp(-\chi)$ and $(L_{|U}, \bD_h)$ is the restriction of
$(L,\bD_h)$ to $U$. Hence $L=\cO_X$ and $RH^{-1}(1)=(\cO_X,0)$.

\medskip

The map $RH$ is surjective, and thus an isomorphism. Indeed, let $\rho$ be a unitary character of $H_1(U)$ corresponding to
 a line bundle with a flat connection $(L_U,\bD)$ on $U$. We define $(L,\alpha)$
as follows. Let $\alpha_i\in [0,1)$ be such that the monodromy of
$(L_U,\bD)$ around a general point of $D_i$ is multiplication by
$\exp(2\pi i\alpha_i)$. Let $(L_U^\vee,\bD^\vee)$ denote the dual connection of $(L_U,\bD)$. Let $(L^\vee,\bD^\vee)$ denote the canonical Deligne extension of $(L_U^\vee,\bD^\vee)$ to a line bundle with a flat connection with logarithmic singularities along $D$.  Then $L$ is the dual of $L^\vee$.  By Proposition \ref{chern}, $c_1(L^\vee)=-\sum_{i\in S}\alpha_i [D_i]$ in $H^2(X,\bR)$. Hence $(L,\alpha)$ is a realization of the boundary $\alpha\cdot D$ such that $RH(L,\alpha)=(L_U,\bD)$. \ \ $\square$

\medskip

\noindent {\bf Comparison with log resolutions.}  Let $X$ be a nonsingular
complex projective variety, $D=\cup_{i\in S}D_i$ a divisor on $X$
with irreducible components $D_i$, and let $U=X-D$. Fix a log
resolution $\mu:Z\ra X$ of $(X,D)$ which is an isomorphism above
$U$. Let $E=Z-U$ with irreducible decomposition $E=\cup_{j\in
S'}E_j$. Let $\mu^*D_i=\sum_{j\in S'}e_{ij}\cdot E_j$. Let $\cV$ be a unitary local system on $U$ and denote by $T_j$ the monodromy of $\cV$ around $E_j$ given by a small loop in $U$ centered at a general point of $E_j$. Since $\cV$ is unitary, $T_j$ is multiplication by $\exp(-i2\pi\beta_j)$ for some $\beta_j\in[0,1)$. If $E_j$ is the
strict transform of some $D_i$, we denote by $T(i)$ the monodromy
$T_j$ and we let $\alpha_i=\beta_j$. Let $\gamma_j$ and $\delta_i$
be the images in $H_1(U)$ of a small loop in $U$ centered at a
general point of $E_j$ and, respectively, $D_i$.
\begin{lem}\label{small}  With the notation as above, for all $j\in S'$,

(a) $\gamma_j=\sum_{i\in S}e_{ij}\delta_i \in H_1(U)$.

(b) $T_j=\prod_{i\in S}T(i)^{e_{ij}}$.

(c) $\beta_j=\{\sum_{i\in S}e_{ij}\alpha_i\}$, the fractional
part.
\end{lem}
\begin{proof} By definitions, (a) implies (b) which implies (c). We
show (a). By duality (e.g. \cite{Bre}-VI. Theorem 8.3), we have
canonical isomorphisms
\begin{align*}
H^{2n-2}(D)= H_2(X,U)\\
H^{2n-2}(E)= H_2(Z,U)
\end{align*}
where $n$ is the dimension of $X$. Since the singular locus
$D_{sing}$ of $D$ has codimension at least 2 in $X$,
$$H^{2n-2}(D)= H^{2n-2}(D,D_{sing})= H_0(D-D_{sing})= \bZ^S.$$
Indeed, $D-D_{sing}$ has one path-connected component for each
$D_i$. Call $\eta_i$ this generator of $H^{2n-2}(D)$, so that this
last space is $\oplus_{i\in S}\bZ\eta_i$. Similarly,
$H^{2n-2}(E)\cong\oplus_{j\in S'}\bZ\xi _j,$ where $\xi_j$
corresponds to $E_j$. The images of $\eta_i$ and $\xi_j$ in
$H_1(U)$ by the boundary maps $H_2(X,U)\ra H_1(U)$ and,
respectively $H_2(Z,U)\ra H_1(U)$, are $\delta_i$ and $\gamma_j$,
respectively. Hence (a) follows if we show that, under the natural
direct image $\mu_*:H_2(Z,U)\ra H_2(X,U)$, there is an equality
$\mu_*(\xi_j)=\sum_{i\in S}e_{ij}\eta_i$.

\medskip

Since $H_2(X,U)$ and $H_2(Z,U)$ are free abelian groups, by the
Universal Coefficient Theorem (\cite{Bre}- V. Corollary 7.2) the
map $\mu_*$ is the adjoint of the pullback map $\mu^*: H^2(X,U)\ra
H^2(Z,U)$. Let $Z_D^{\; 2}(X)$ and $Z^{\; 2}_E(Z)$ be the groups
of $2$-codimensional cycles on $X$ and $Z$ with supports in $D$
and $E$, respectively. There are cycle maps (\cite{Fu} -p. 380)
$cl^X$ and $cl^Z$ making the following diagram commutative
(\cite{Fu} -Proposition 19.2; see also Example 19.2.6):
$$
\xymatrix{
Z_D^{\; 2}(X)  \ar[r]^{cl^X} \ar[d]_{\mu^*}&   H^2(X,U) \ar[d]_{\mu^*}  \\
Z^{\; 2}_E(Z)  \ar[r]^{cl^Z} &  H^2(Z,U)\\
}
$$
By definition in loc. cit., $cl^X(D_i)=\eta_i^\vee$, under the
isomorphism $$H^2(X,U)=H_{2n-2}(D)=\oplus _iH_{2n-2}(D_i)=\oplus
_iH^{2n-2}(D_i)^\vee.$$ Similarly, $cl^Z(E_j)=\xi_j^\vee$. Since
the cycle pullback $\mu^*D_i$ is $\sum_{j\in S'}e_{ij}E_j$, the
claim follows.
 \ \
\end{proof}

The canonical generators $\eta_i$ and $\eta_i^\vee$ of $H_2(X,U)$
and $H^2(X,U)$ will be called $[D_i]^\vee$ and $[D_i]$,
respectively, when no confusion arises.

\begin{lem}\label{de} With the notation as in Lemma \ref{small}, let $(M_U,\bD)$ be a flat connection on $U$ corresponding to the local system $\cV$. Denote by $(M,\bD)$ the canonical Deligne extension to $Z$. There exists a line bundle $L$ on $X$ such that the
canonical Deligne extension is
$${M}=\mu^*L\otimes\cO_Z(\rndown{e}\cdot E),$$
where $e=(e_j)_{j\in S'}$ with $e_j=\sum_{i\in S}e_{ij}\alpha_i$.
\end{lem}
\begin{proof} Recall the local description of ${M}$. Let $z_j$ be
local coordinates of $E_j$. Around a general point on $E_j$, let
$u_j$ be a frame of the multi-valued flat sections of
$(M_U,\bD)$. Then $v_j=\exp(\log z_j\cdot\beta_j)\cdot u_j$ is a
single-valued holomorphic section of $M_U$. Indeed,
$$T_jv_j=\exp((\log z_j+2\pi i)\beta_j)\cdot \exp(-2\pi
i\beta_j)u_j=v_j.$$ Then ${M}$ is the line bundle obtained by
declaring $v_j$ to be a local holomorphic frame on $Z$ and gluing.
Define
$$w_j=\exp(\log z_j\cdot e_j)\cdot u_j,$$
so that by Lemma \ref{small}-(c), $v_j=\exp(\log z_j\cdot
(-\rndown{e_j}))\cdot w_j$. Then
$$T_jw_j=\exp(2\pi i e_j)\cdot\exp(-2\pi i\beta_j)\cdot
w_j=w_j.$$ Hence $w_j$ is also single-valued. Moreover,  $w_j$ is
the restriction to an appropriate open subset of $Z$ of
$w'\circ\mu$, where $w'=\exp(\sum_{i\in S}\log f_i(y)\cdot
\alpha_i)\cdot u'$ with $f_i$ a local equation on $X$ for $D_i$
and $u'$ a local frame for the multi-valued flat sections of
$(M,\bD)$ with respect to a chart (with coordinates $y$) in
$X$. Hence $w'$ defines, locally in charts on $X$, a holomorphic
frame for $M_U$ on $U$. We let $L$ to be  the line bundle
on $X$ obtained by declaring a local frame to be $w'$ and then
gluing. The Lemma follows.\ \
\end{proof}

\begin{prop}\label{comparison}
With the notation as in Lemma \ref{small}, the map $$\mu^*_{par}:\Pic^\tau(X,D) \ra\Pic^\tau(Z,E)$$ given by
$(L,\alpha)  \mapsto (\mu^*L-\rndown{e}\cdot E, \{e\})$ is an isomorphism, where $e\in \bR^{S'}$ is given by $\mu^*(\alpha\cdot
D)=e\cdot E$.
\end{prop}
\begin{proof} It is clearly an injective group homomorphism. To show that $\mu^*_{par}$ is surjective, let $(M^{-1},\beta)\in\Pic^\tau(Z,E)$ for some line bundle $M$ on $Z$. Let $\alpha$ be defined in terms of $\beta$ as in Lemma \ref{small}. By the proof of Theorem \ref{ls} for the simple normal crossings case,  $(M^{-1},\beta)$ corresponds under the isomorphism $RH$ to the flat connection $(M^{-1}_{|U},\bD)$ such that $M$ is the  canonical Deligne extension of the dual connection. By Lemma \ref{de}, $M=\mu^*L\otimes\cO_Z(\rndown{e}\cdot E)$ for some line bundle $L$ on $X$.
By Lemma \ref{small}-(c),  $(L^{-1},\alpha)$ lies in $\Pic^\tau(X,D)$. Thus $\mu^*_{par}(L^{-1},\alpha)=(M^{-1},\beta)$.
\end{proof}

\medskip
\noindent {\bf Proof of Theorem \ref{ls} for the general case.} It follows by the simple normal crossings case plus Proposition \ref{comparison}.\ \ $\square$

\section{Naturality and structure}\label{ns}

\noindent {\bf Naturality.} We show now that the canonical
isomorphism RH is natural in the sense which we describe
below.  Let $X$ be a nonsingular complex projective variety,
$D=\cup_{i\in S}D_i$ a divisor on $X$ with irreducible components
$D_i$, and let $U=X-D$. Consider the long exact sequence of
homology with integral coefficients
\begin{equation}\label{leh} \cdots\ra
H_2(X)\xrightarrow{\psi} H_2(X,U)\ra H_1(U)\xrightarrow{\phi}
H_1(X) \ra 0
\end{equation}
where the maps $\psi$ and $\phi$ are defined as illustrated.
Because $S^1$ is injective, there is a short exact sequence
\begin{align}\label{comp1}
0\ra\Hom (H_1(X),S^1)\ra\Hom (H_1(U),S^1)\ra \Hom
\left({{\rm{ker}}(\phi)}, S^1\right)\ra 0.
\end{align}
On the other hand, there is a short exact sequence
\begin{align}\label{comp2}
0\ra \Pic^\tau(X)\ra\Pic^\tau(X,D)\ra B(X,D)\ra 0,
\end{align}
where the injective morphism is $L\mapsto (L,0)$.

\begin{prop}\label{comp} With the notation as above, the  sequences
(\ref{comp1}) and (\ref{comp2}) are compatible via the isomorphism
RH.
\end{prop}
\begin{proof}
 Let $H_2(X,U)=\oplus_{i\in S}\bZ
[D_i]^\vee$ as in the proof of Lemma \ref{small}. Let $K$ be the
group $\Hom\left ({{\rm{ker}}(\phi)},S^1\right )$. Then $K$ is the
kernel of the map $\Hom(H_2(X,U), S^1)\ra \Hom ({\rm{Im}} (\psi)
,S^1)$, hence
\begin{align*}
K=\{ \alpha\in [0,1)^S :  \sum_{i\in S}\beta_i\alpha_i \in \bZ\
{\rm{whenever\ }} \sum_{i\in S}\beta_i[D_i]^\vee \in
{\rm{Im}}(\psi),\ \beta_i\in\bZ\ \}.\\
\end{align*}
Consider the map $\Pic^\tau(X,D)\ra K$ given by $(L,\alpha)\mapsto
\alpha$. By Lemma \ref{wd} below $B(X,D)\subset K$, i.e. this map
is well-defined. The kernel is  $\Pic^\tau(X)$. This gives a
morphism (i.e. commuting diagram) from the exact sequence
(\ref{comp2}) to (\ref{comp1}), where the first two vertical maps
are RH (the first map is for the case $D=0$), and the last
vertical map is the inclusion $B(X,D)\subset K$. By the snake
lemma, $B(X,D)\cong K$ canonically, i.e. we obtain a short exact
sequence
\begin{equation}\label{ker}
0\ra \Pic^\tau(X)\ra\Pic^\tau(X,D)\ra K\ra 0
\end{equation}
\ \
\end{proof}

\begin{lem}\label{wd} With the notation as in Proposition \ref{comp},
let $[D_i]$ be the real cohomology classes of the divisors $D_i$,
$[D_i]^\vee$ the canonical generators of the integral homology
group $H_2(X,U)$. Let $a=\sum_{i\in S}\alpha_i [D_i]$ be an
element in the image of $H^2(X,\bZ)$ in $H^2(X,\bR)$, with
$\alpha_i\in\bR$ (hence, actually $\alpha_i\in\bQ$). Let $b=\sum_{i\in S}\beta_i[D_i]^\vee$ be an
element in the image of $\psi$, with $\beta_i\in\bZ$. Then
$\sum_{i\in S}\beta_i\alpha_i \in \bZ$ and depends only on $a$ and
$b$.
\end{lem}
\begin{proof}
This is a corollary of the existence of the cap product. Consider
the cap product for integral coefficients
$$H^2(X)\otimes H_2(X)\xrightarrow{\cap} H_0(X)\cong\bZ.$$
If $t\in H^2(X)$ is a torsion class, then $mt=0$ for some
$m\in\bZ_{>0}$. Hence $m(t\cap b)$ and $t\cap b$ are $0$ for all
$b\in H_2(X)$. The group $H^2(X)$ modulo torsion can be identified
with the image of $H^2(X)\ra H^2(X,\bR)$. Let $A$ be the free
finitely generated subgroup of $H^2(X)/{\rm{torsion}}$ given by
the intersection with the subspace $V$ generated by the $[D_i]
(i\in S)$ in $H^2(X,\bR)$. Then the cap product induces a bilinear
map
$$A\otimes H_2(X)\ra\bZ.$$
If $u$ is in the image of $H_2(U)\ra H_2(X)$ and $a\in A$, then
$a\cap u =0$. Hence we have a map
\begin{align}\label{pairing}
A\otimes B\ra \bZ,
\end{align}
where $B$ is the quotient of $H_2(X)$ by the image of $H_2(U)$.
From (\ref{comp1}), $B$ can be identified with the image of
$\psi$. This proves the statement of the Lemma for the case when
$\alpha_i$ are integers.

\medskip

It remains to show that the bilinear map (\ref{pairing}) extends
to a bilinear map
$$l^{-1}(A)\otimes B\ra \bZ,$$
where $l:\bR^S\ra V$ is the linear map
$\alpha\mapsto\alpha\cdot[D]$. Let $W$ be the subspace $l^{-1}(0)$
of $\bR^S$. Then it is enough to show that the pairing induced by
linearity between $W$ and $B$ is zero. By (\ref{pairing}) this is
true for the subgroup $W\cap \bZ^S$ in $W$. The claim follows
since $W\cap \bZ^S$ is a lattice in $W$.\ \
\end{proof}

\begin{lem}\label{torsion} With the notation as in Proposition
\ref{comp}, let $N>1$ be an integer. There is a natural
isomorphism between the group of characters of $H_1(U,\bZ_N)$ and
the subgroup of $N$-torsion elements of $\Pic^\tau(X,D)$.
\end{lem}
\begin{proof}
By RH, the subgroup of $N$-torsion elements of
$\Pic^\tau(X,D)$ can be identified with $\Hom (H_1(U),\mu_N)$,
where $\mu_N=\Hom (\bZ_N,S^1)$. By adjunction of functors, the last group is
naturally isomorphic with $\Hom (H_1(U)\otimes_\bZ\bZ_N,S^1)$.
\end{proof}

\medskip

\noindent{\bf The realizable boundaries.} We consider the set of
realizable boundaries $B(X,D)$. By the previous subsection,
$B(X,D)$ is the direct sum of a finite group with $(S^1)^{r(U)}$,
where $r(U):=b_1(U)-b_1(X)$. Indeed, $B(X,D)$ is canonically
isomorphic with $\Hom (\ker (\phi),S^1)$ by Proposition
\ref{comp}, and $\ker(\phi)$ is a finitely generated abelian group
of rank $r(U)$. Remark that the first Betti number of $X$,
$b_1(X)$, is  a birational invariant since it equals twice
the Hodge number $h^{1,0}(X)$. By definition, $B(X,D)$ has an
inclusion in $\bR^S\cap [0,1)^S$.

\begin{lem}\label{rebd} With the notation as in Proposition
\ref{comp}, $B(X,D)$ is a finite disjoint union of rational convex
polytopes in $\bR^S$. More precisely, for the linear map $l$ of
$\bR^S$ onto $\bR^{|S|-r(U)}$ defined in the proof of Lemma
\ref{wd}, there is a finite set of points $p_k\in\bR^{|S|-r(U)}$
such that $B(X,D)$ is the union of the polytopes
$P_k=l^{-1}(p_k)\cap[0,1)^S$.
\end{lem}
\begin{proof}
With the notation of the proof of Lemma \ref{wd},
$$B(X,D)=l^{-1}(\Lambda)\cap [0,1)^S,$$
where $\Lambda$ is the subset of $V$ consisting only of cohomology
classes which are realized as first Chern classes of line bundles
on $X$. Then $\Lambda$ is a free finitely generated abelian group
and $\Lambda\otimes_\bZ\bR=V$. That is $\Lambda$ is a lattice in
$V$. There are only finitely many lattice points in $V$ in the
image of $[0,1)^S$ under $l$.
\end{proof}

\medskip

\noindent{\bf Structure.} It follows from the short exact sequence
(\ref{ker}) and the construction of the map $l$ of Lemma \ref{wd}
that, topologically, $\Pic^\tau(X,D)$ has a finite disjoint
decomposition
\begin{align}\label{desc}
\Pic^\tau(X,D)=\coprod_k P_k\times \Pic^\tau(X),
\end{align}
where $P_k$ are rational convex polytopes  in $\bR^S\cap [0,1)^S$.
Pointwise, the subset of $\Pic^\tau(X,D)$ corresponding to
$P_k\times \Pic^\tau(X)$ consists of realizations $(L_k+M,\alpha)$
of boundaries $\alpha\in P_k$, where $M\in\Pic^\tau(X)$, for some
line bundle $L_k$.

\begin{prop}\label{struc} For each $k$, $L_k$ can be chosen such
that $(L_k,\alpha_k)$ is torsion, for some $\alpha_k\in P_k$.
\end{prop}
\begin{proof} Let $\alpha\in P_k$ be a rational point. It is enough
to show that there exists some $(L,\alpha)\in \Pic^\tau(X,D)$
which is torsion. Fix $L$ such that $(L,\alpha)\in\Pic^\tau(X,D)$.
It is enough to show that $(L,\alpha)\cdot\Pic^\tau(X)$ contains a
torsion element.

From the exponential exact sequence we have the short exact
sequence
$$0\ra\Pic^0(X)\ra\Pic^\tau(X)\ra H^2(X,\bZ)_{tor}\ra 0.$$
Hence we can write
$$\Pic^\tau(X)=\coprod_{E} E\cdot\Pic^0(X),$$
where the disjoint union is taken over a fixed set of
representatives for $H^2(X,\bZ)_{tor}$ in $\Pic^\tau(X)$. Thus
$$(L,\alpha)\cdot\Pic^\tau(X)=\coprod_{E}(L+E,\alpha)\cdot\Pic^0(X).$$
It is enough for our purpose to show that $(L+E,\alpha)\cdot
\Pic^0(X)$ contains a torsion element of $\Pic^\tau(X,D)$.

\medskip

Let $g=(L+E,\alpha)$ and $H=\Pic^0(X)$. Since $\alpha\in\bQ^S$ and
$H^2(X,\bZ)_{tor}$ is finite, there exists $N>1$ such that $g^N\in
H$. Since $H$ is divisible, there exists $h\in H$ such that
$h^N=(g^{-1})^N$. Thus $(gh)^N=1$ and $gh$ is torsion.$\ \
 $
\end{proof}

\begin{example}\label{ex1} Consider the case of $(\bP, D)$, where $\bP=\bP^2$ and $D$ is
the union of 5 distinct lines $D_i$ ($i=1, .., 5$) intersecting in
one point in $\bP$. Here $\Pic^\tau(\bP)=\{\cO\}$ and
$\Pic^\tau(\bP,D)=B(\bP,D)$. The map $l$ from Lemma \ref{rebd} is
$l(\alpha)=(\sum_{i=1,..,5}\alpha_i)[D_1]$, $p_k=k[D_1]$ . The
canonical decomposition of $B(\bP,D)$ is $\coprod_{k=0, .., 4}
P_k$, where $P_k=[0,1)^5\cap l^{-1}(p_k)$ with $k=0,..,4$.
\end{example}

\section{Finite abelian covers}\label{facs}

 The subject of finite
abelian covers is well-known and has been studied by many people.
To our knowledge, the geometric characterization of abelian covers
has two versions which appear in their final form in \cite{Na} and
\cite{Pa}. The isomorphism of Theorem \ref{ls} recovers both
characterizations and allows us to prove some results which, to
our knowledge, appear in their published version only under some
additional hypotheses.

\medskip

Let $X$ be a nonsingular complex projective variety, $D=\cup_{i\in
S}D_i$ a divisor on $X$ with irreducible components $D_i$, and let
$U=X-D$. Let $G$ be a finite abelian group.

\begin{df}\label{gcovering} A map $Y\ra X$ is a {\it $G$-cover} if it is a finite map
together with a faithful action of $G$ on $Y$ such that the map
exhibits $X$ as the quotient of $Y$ via $G$. Two covers are
said two be equivalent
 if there is an isomorphism between them commuting
with the cover maps.
\end{df}

Recall the following topological characterization (see Appendix 1
to Chapter VIII of \cite{Za}, or \cite{Na}). The morphisms of
$H_1(U)$ onto $G$ are in one-to-one correspondence with the
equivalence classes of unramified $G$-covers of $U$. These, in
turn, are in one-to-one correspondence with equivalence classes of
normal $G$-covers of $X$ unramified above $U$. The group $G$ is
recovered as the group of automorphisms of the cover commuting
with the cover map.

\medskip

\noindent {\bf Proof of Corollary \ref{fc}.} By the topological
characterization and by the canonical isomorphism RH, it is enough
to show there is a one-to-one equivalence between surjections
$H_1(U)\ra G$ and subgroups $G^*\subset \Hom(H_1(U),S^1)$. This is
a standard exercise in duality.
 \ \
$\square $


\begin{rem}\label{ggg} It is well-known that the local systems corresponding
to $G^*\subset \Hom(H_1(U),S^1)$ in Corollary \ref{fc} can also be
obtained as follows. Let $\pi_U:V\ra U$ be the unramified cover
induced by a surjection $H_1(U)\ra G$. Then $(\pi_U)_*\bC_V$ is a
higher rank local system on $U$ which decomposes into rank one
unitary local systems $\cV_\chi$ $(\chi\in G^*)$. Moreover, there
is an eigensheaf decomposition of $(\pi_U)_*\cO_V$ into line
bundles $\cM_\chi$ ($\chi\in G^*$), and a natural derivative
$\nabla$ such that $(\cM_\chi,\nabla)$ is the flat connection
corresponding to $\cV_\chi$.
\end{rem}


\medskip

\noindent{\bf Recovering the geometric characterization of
\cite{Pa}.} Let $X, D, U$, and $G$ be as above. Pardini \cite{Pa}
gives the following geometric characterization of finite abelian
covers which also follows from Theorem \ref{ls}.

\begin{df} The data $\{ (L_\chi, H_i, \psi_i) : \chi\in G^*, i\in S
\}$  is called {\it building data on $X$ with branch $D$} if $H_i$
is a cyclic subgroup of $G$, $\psi$ is a generator of $H_i^*$,
$L_\chi$ are line bundles satisfying the linear relation
\begin{align}\label{lineq}
L_\chi+L_{\chi'}=L_{\chi\chi'}+\varepsilon_{\chi,\chi'}\cdot D,
\end{align}
where $\varepsilon_{\chi,\chi'}\in \{0,1\}^S$ with
$$(\varepsilon_{\chi,\chi'})_i=
\begin{cases}
0 & \text{if }\iota_{\chi,i}+\iota_{\chi',i}< m_i,\\
1 & \text{otherwise},
\end{cases}
$$
Here $m_i$ is the order of $H_i$, and
$\iota_{\chi,i}\in\{0,\ldots,m_i-1\}$ is given by
$\chi_{|H_i}=\psi_i^{\iota_{\chi,i}}$.
\end{df}

\begin{cor}\label{pardini} There is a one-to-one
correspondence between the set of equivalence classes of normal
$G$-covers of $X$ unramified above $U$ and building data on $X$
with branch $D$. This correspondence is the same as the one of
\cite{Pa}.
\end{cor}

Pardini's correspondence is as follows. Suppose we start with a
normal $G$-cover $\pi:Y\ra X$ unramified above $D$. Then there
is an eigensheaf decomposition
$$\pi_*\cO_Y=\bigoplus_{\chi\in G^*}L_\chi^{-1}.$$
Let $H_i$ be the inertia group of (any=all) components $T$ of
$\pi^{-1}(D)$. Let $\psi_i$ be the representation of $H_i$ on
$\bf{m}/\bf{m}^2$ induced by the cotangent map, where $\bf{m}$ is
the maximal ideal of the local ring $\cO_{Y,T}$. Then $(L_\chi,
H_i,\psi_i)$ is the corresponding building data. Conversely, suppose we start with building data $(L_\chi,
H_i,\psi_i)$. Define the $\cO_X$-linear multiplication maps
$\mu_{\chi,\chi'}:L_\chi^{-1}\otimes L_{\chi'}^{-1}\ra
L_{\chi\chi'}^{-1}$ by setting
$$\mu_{\chi,\chi'}=\prod_{i\in
S}\sigma_i^{(\varepsilon_{\chi,\chi'})_i},$$ viewed as a global
section of $L_\chi\otimes L_{\chi'}\otimes L_{\chi\chi'}^{-1}$,
where $\sigma_i$ are sections of $\cO_X(D_i)$ vanishing on $D_i$.
The corresponding $G$-cover is then defined (up to equivalence
due to freedom of rearranging the $\chi$'s around) as
\begin{align}\label{decomp}
Y=\text{\bf Spec }_{\cO_X}(\bigoplus_{\chi\in G^*}L_\chi^{-1}).
\end{align}

{\it Proof of Corollary \ref{pardini}.} By Corollary \ref{fc}, the
 normal $G$-covers of $X$ unramified
above $U$ up to equivalence are into one-to-one correspondence
with subgroups $G^*\subset\Pic^\tau(X,D)$. Start with
$G^*=\{(L_\chi,\alpha_\chi) : \chi\in G^*
\}\subset\Pic^\tau(X,D)$, corresponding to the $G$-cover
$\pi:Y\ra X$ and to the epimorphism $\rho:H_1(U)\ra G$. To this we
attach building data $(L_\chi, H_i,\psi_i)$ as follows.

\medskip

Let $\delta_i$ be homology class of a small loop in $U$ centered
at a general point of $D_i$. Let $H_i$ be the cyclic subgroup of
$G$ generated by $\rho(\delta_i)$. Let $\psi_i$ be the character
of $H_i$ taking $\rho(\delta_i)$ to $\exp(i2\pi(1/m_i))$, where
$m_i$ is the order of $H_i$. Then $(L_\chi,H_i,\psi_i)$ is a
building data. Indeed,  $(L_\chi,\alpha_\chi)$ corresponds to the
local system given by $\chi\circ\rho$. Hence, by the
description of the isomorphism RH, $\chi$ sends $\rho(\delta_i)$ to
$\exp(i2\pi\alpha_{\chi,i})$. Thus
$\iota_{\chi,i}=m_i\alpha_{\chi,i}$ and the linear relation
(\ref{lineq}) follows from the group operation on
$\Pic^\tau(X,D)$.

\medskip

Let $\pi':Y'\ra X$ be the $G$-cover corresponding to the
building data $(L_\chi, H_i,\psi_i)$  via Pardini's
correspondence. It is enough to show that  $\pi'$ is equivalent to
$\pi$. It is actually enough to show that the corresponding
unramified covers above $U$ are equivalent. This follows from
Remark \ref{ggg} and the explicit description of the isomorphism
RH . \ \  $\square$

\smallskip

\noindent {\bf Proof of Corollary \ref{ss}.} It follows from the Corollary \ref{pardini} and (\ref{decomp}).\ \ $\square$

\smallskip

\noindent{\bf Proof of Corollary \ref{ss2}.} Let $E=Z-U$ and $E_j\
(j\in S')$ be its irreducible components. We have obtained in
Proposition \ref{comparison} that the map $$\mu^*_{par}:\Pic^\tau(X,D)  \ra\Pic^\tau(Z,E)$$
given by $(L,\alpha)  \mapsto (\mu^*L-\rndown{e}\cdot E, \{e\})$
is an isomorphism, where $e\in \bR^{S'}$ is given by $\mu^*(\alpha\cdot
D)=e\cdot E$. Therefore Corollary \ref{ss2} follows from
Corollary \ref{fc} and Corollary \ref{ss}.\ \ $\square$

\smallskip

\noindent {\bf Proof of Corollary \ref{ss3}.} The space of global
$q$-forms is a birational invariant of nonsingular projective
varieties. Hence we can choose any desingularization $Y$ of
$\wti{X}$. Let $Y$ and maps $\eta:Y\ra \wti{X}$, $\nu:Y\ra\wti{Z}$
be a common $G$-equivariant desingularization (\cite{RY},
\cite{AW}) of $\wti{X}$ and $\wti{Z}$, so that we have a
commutative diagram
$$
\xymatrix{
& Y \ar[dl]_\nu \ar[dr]^\eta \ar@/_/[ddr]_f & \\
\wti{Z} \ar[d]_\rho & & \wti{X} \ar[d]^\pi \\
Z \ar[rr] ^\mu & & X }$$ where $f$ is the induced map from $Y$ to
$X$. Corollary \ref{ss3} follows from Proposition \ref{equiv} for
the case when $M=\cO_X$ and by Hodge duality. \ \ $\square$

\begin{prop}\label{equiv}
With the notation as in the proof of Corollary \ref{ss3}, let $M$
be a line bundle on $X$. Let $\chi\in G^*$ correspond to
$(L_\chi,\alpha_\chi)\in\Pic^\tau(X,D)$. Then
$$H^q(Y,f^*M^{-1})_\chi\cong H^{n-q}(X,\omega_X\otimes M\otimes
L_\chi\otimes \cJ(\alpha_\chi\cdot D)),$$ where the subscript
$\chi$ denotes the $\chi$-eigenspace and $n$ is the dimension of
$X$.
\end{prop}
\begin{proof} By Lemma 3.24 of \cite{EV2}, $\wti{Z}$ has quotient singularities. Hence, by Proposition 1 of \cite{Vie}, $\wti{Z}$ has rational singularities.
That is $\nu_*\cO_Y=\cO_{\wti{Z}}$ and $R^i\nu_*\cO_Y=0$ for
$i>0$. Then, by the projection formula (\cite{H}-Ex. III.8.3),
$\nu_*(f^*M^{-1})=\rho^*\mu^*M^{-1}$ and $R^i\nu_*(f^*M^{-1})=0$
for $i>0$. Hence, by the equivariant version of the Leray spectral
sequence (\cite{La}-Proposition B.1.1), there is a $G$-equivariant
isomorphism
$$H^q(Y,f^*M^{-1})\cong H^q(\wti{Z},\rho^*\mu^*M^{-1})$$ for all $q$.
Now, $\rho$ is a finite morphism, so $R^i\rho_*\cO_{\wti{Z}}=0$
for $i>0$. By the projection formula again,
$R^i\rho_*(\rho^*\mu^*M^{-1})$ vanishes for $i>0$ and is equal to
$\rho_*\cO_{\wti{Z}}\otimes\mu^*M^{-1}$ for $i=0$. Then the
equivariant version of \cite{La}-Proposition B.1.1 gives again a
$G$-equivariant isomorphism
$$H^q(\wti{Z},\rho^*\mu^*M^{-1})\cong H^q(Z,\rho_*\cO_{\wti{Z}}\otimes\mu^*M^{-1})$$ for all $q$. Hence, by
Corollary \ref{ss2}, we have an isomorphism
$$H^q(\wti{Z},\rho^*\mu^*M^{-1})_\chi\cong H^q(Z,\mu^*(M^{-1}\otimes L_\chi ^{-1})\otimes\cO_Z(\rndown{\mu^*(\alpha_\chi\cdot D)}))=\star.$$
By Serre duality on $Z$,
$$\star\cong H^{n-q}(Z,\omega_Z\otimes\mu^*(M\otimes L_\chi)\otimes\cO_Z (-\rndown{\mu^*(\alpha_\chi\cdot D)})).$$
By Theorem \ref{localvanishing}, we are again in the situation of
\cite{La}-Proposition B.1.1 and obtain
$$\star\cong H^{n-q}(X,\omega_X\otimes M\otimes
L_\chi\otimes\cJ(\alpha_\chi\cdot D)).\ \  $$
\end{proof}

\section{Canonical stratifications}\label{canstrat}

To prove Theorem \ref{gv} we use the following version due to C.
Simpson \cite{Si} which builds on a previous result of
Green-Lazarsfeld \cite{GL2} and D. Arapura \cite{Ar1}.

\begin{thm}\label{si} (Simpson)  Let $f:Y\ra X$ be a morphism between two
nonsingular complex projective varieties.  Let $G$ be a finite
abelian group acting by automorphisms on $Y$, trivially on $X$,
and suppose $f$ is $G$-equivariant. Fix $\chi\in G^*$ and integers
$q$, $i$. Then the set
$$V^q_i(f,\chi):=\{ M\in \Pic^\tau(X) : h^q(Y,f^*M)_\chi \ge i \}$$
is a finite union of torsion translates of complex subtori of
$\Pic^\tau(X)$, and so is any intersection of these translates.
\end{thm}

This is the "Higgs field equals zero" case (see \cite{Si}-section
5) of a stronger equivariant version which is not explicitly
stated in \cite{Si} but which follows straight-forwardly from the
results there. Indeed, one only needs the equivariant version of
\cite{Si}-Proposition 7.9 to construct the appropriate absolute
functor (in the terminology of loc. cit.) and, thus, the
appropriate absolute closed subset of the moduli space of rank one
local systems on $X$. By the definition in \cite{Si}-section 6,
intersections of absolute closed subsets are again absolute closed
subsets. Then \cite{Si}-Theorem 6.1 applies and the restriction to
$\Pic^\tau(X)$ gives Theorem \ref{si} as stated. See also Theorem
\ref{sigen}. For sake of completeness, we give a less conceptual proof of Theorem \ref{si} based on the trick of \cite{Ar1} (see also \cite{Si}- section 5) where we replace Hodge decomposition by eigenspace decomposition.

\medskip

\noindent{\it Proof of Theorem \ref{si}.}  Let $V^q_i(f)=\{ M\in
\Pic^\tau(X)\ :\ h^q(Y,f^*M)\ge i\ \}$. By \cite{Si}- section 5,
the conclusion of the theorem holds for $V^q_i(f)$. Since
$h^q(Y,f^*M)=\sum_{\chi\in G^*}h^q(Y,f^*M)_\chi$, one has the
equality
$$
V^q_i(f)=\mathop{\bigcup_{P\; :\; G^*\ra \bN}}_{\sum_\chi P(\chi)=i} \left [\bigcap_{\chi\in G^*} V^q_{P(\chi)}(f,\chi)\right ].
$$
Let $V$ be an irreducible component of the closed subvariety
$V^q_i(f,\chi)$ of $\Pic^\tau(X)$. To finish the proof, it is
enough to show that $V$ is an irreducible component of a set of
the type $V^q_I(f)$ for some $I$. Define a function $P :
G^*\ra\bN$ , dependent on $q$, by $P(\psi)=\max\{ j\ |\ V\subset
V^q_j(f,\psi)\}$. Then $V$ is an irreducible component of
$\cap_{\psi\in G^*}V^q_{P(\psi)}(f,\psi)$ since
$V^q_{P(\chi)}(f,\chi)\subset V^q_i(f,\chi)$. But $V$ is not
included in $\cap_{\psi\in G^*} V^q_{P'(\psi)}(f,\psi)$ for any
other function $P': G^* \ra \bN$ with $\sum_{\psi}P
'(\psi)=\sum_\psi P (\psi)=:I$. Hence $V$ is an irreducible
component of $V^q_I(f)$.\ \ $\square$

\begin{lem}\label{si2} Let $X$ be a nonsingular complex projective
variety, $D=\cup_{i\in S}D_i$ a divisor on $X$ with irreducible
components $D_i$. Let $(L,\alpha)$ be a torsion element in
$\Pic^\tau(X,D)$. Then the set
$$V_i^q(L,\alpha):=\{ M\in \Pic^\tau(X) : h^q(X,\omega_X\otimes
M\otimes L\otimes \cJ(\alpha\cdot D))\ge i \}$$ is a finite union
of torsion translates of complex subtori of $\Pic^\tau(X)$, and so
is any intersection of these translates.
\end{lem}
\begin{proof} Choose any finite group $G$ with an embedding
$G^*\subset \Pic^\tau(X,D)$ (hence $G$ is abelian) such that $(L,\alpha)\in G^*$. Let
$\chi\in G^*$ be the character corresponding to $(L,\alpha)$ and
construct $f$ as in Proposition \ref{equiv}. By Proposition
\ref{equiv},
$$V_i^q(L,\alpha)=\{ M^{-1} : M\in V_i^{n-q}(f,\chi) \}.$$
The Lemma follows then from Theorem \ref{si}.
\end{proof}

\smallskip

Let $X$ be a nonsingular complex projective variety of dimension
$n$, $D=\cup_{i\in S}D_i$ a divisor on $X$ with irreducible
components $D_i$, and let $U=X-D$. Let $G$ be a finite abelian
group. Define
$$V^q_i(X,D):=\{ (L,\alpha)\in \Pic^\tau(X,D) : h^q(X, \omega_X\otimes
L\otimes\cJ (\alpha\cdot D))\ge i \},$$ where $\cJ (\alpha\cdot
D))$ denotes the multiplier ideal of the $\bR$-divisor
$\alpha\cdot D$. Whenever the context leaves no room for
ambiguity, we will just write $V^q_i$ instead of $V^q_i(X,D)$.

\begin{lem}\label{indep} With the notation as above, the set
$V^q_i(U)$ of unitary rank one local systems on $U$ corresponding
to $V^q_i(X,D)$ under the isomorphism RH depends only on
$U$ and not on $(X,D)$.
\end{lem}
\begin{proof} Let $\mu:(Z,E)\ra (X,D)$ be a
log-resolution which is an isomorphism above $U$, where
$E=\cup_{j\in S'}E_j$ is the inverse image $\mu^{-1}(D)$. The map
$$\mu^*_{par}: \Pic^\tau(X,D)  \ra\Pic^\tau(Z,E)$$
given by $(L,\alpha) \mapsto (\mu^*L-\rndown{e}\cdot E, \{e\})$ is an isomorphism, where $e\in \bR^{S'}$ is given by $\mu^*(\alpha\cdot
D)=e\cdot E$. By Theorem \ref{localvanishing},
$R^j\mu_*(\omega_{Y/X}\otimes\cO_Y(-\rndown{\mu^*(\alpha\cdot
D)}))=0$ for $j>0$, and it follows that under the map $\mu^*_{par}$ the sets
$V^q_i(X,D)$ and $V^q_i(Z,E)$ are into one-to-one equivalence.
Since the map $\mu^*_{par}$ commutes with the isomorphisms
RH, $V^q_i(X,D)$ and $V^q_i(Z,E)$ induce the same subset
of the space of unitary rank one local systems on $U$. The case of
two different divisorial compactifications of $U$ is reduced to
the above case by considering a common log resolution.
\end{proof}

 We can actually describe $V^q_i(U)$ purely
in terms of local systems. For a unitary local system $\cV$ on
$U$, let $F$ be the Hodge filtration on $H^*(U,\cV)$ constructed
by Timmerscheidt \cite{Ti}.

\begin{prop}\label{timm}
With the notation as in Lemma \ref{indep},
$$V^q_i(U)=\{\ \cV\in\Hom (H_1(U),S^1)\ |\ \dim \Gr_F^{0}H^{n-q}(U,\cV^\vee)\ge i\ \}.$$
\end{prop}
\begin{proof} By Lemma \ref{indep}, we may assume $U=X-D$ with $D$ having simple normal crossings. Let $\cV\in \Hom (H_1(U),S^1)$ and let $(L,\alpha)\in \Pic^\tau (X,D)$ be the corresponding realization of boundary. Recall from the proof of Theorem \ref{ls} that $L^{-1}$ is the line bundle of the canonical Deligne extension to $X$ of the connection on $U$ given by the dual local system $\cV^\vee$. We have by \cite{Ti} - 2nd Theorem, part (a),
$$\dim \Gr_F^0 H^{n-q}(U,\cV^\vee)=\dim H^{n-q}(X, \Omega^0_X(\log D)\otimes L^{-1})=\star  ,$$
where $\Omega^p_X(\log D)$ denotes the sheaf of $p$-differential forms with logarithmic poles along $D$, so that $\Omega^0_X(\log D)=\cO_X$ by definition. Then, by Serre duality,
$$\star = \dim H^q(X, \omega_X\otimes L)$$
and this last space is the same as $H^q(X,\omega_X\otimes L\otimes \cJ (\alpha\cdot D))$ since $\cJ(\alpha\cdot D)=\cO_X$ in this case.
\end{proof}

\smallskip

\noindent{\bf Proof of Theorem \ref{gv}.} From (\ref{desc}) and
Proposition \ref{struc}, we can write $V_i^q$  as a disjoint union
\begin{align*}
V_i^q & =\coprod_{k}
V_i^q\cap \left [P_k\times \Pic^\tau(X)\right ]\\
&=\coprod_k V^q_i\cap [\{(L_k,\alpha) :  \alpha\in P_k \}\cdot
\Pic^\tau(X)],
\end{align*}
where $P_k$ are rational convex polytopes in $\bR^S$, and
$(L_k,\alpha_k)$ is torsion in $\Pic^\tau(X,D)$ for some
$\alpha_k\in P_k$.

\medskip

\noindent{\it The simple normal crossings case.} Assume now that
$D$ has simple normal crossings. Then
$V^q_i(L_k,\alpha_k)=V^q_i(L_k,\alpha)$ for all $\alpha\in P_k$
since $\cJ(\alpha\cdot D)=\cO_X$. We will call this set
$V^q_i(L_k)$. Hence,
\begin{align*}
V^q_i  =\coprod_k P_k\times V^q_i(L_k) =\coprod_k \{(L_k,\alpha) : \alpha\in P_k \}\cdot V^q_i(L_k).
\end{align*}
Then Theorem \ref{gv} follows in this case from Lemma \ref{si2}.

\medskip

\noindent{\it General case.} Assume now that $(X,D)$ is as in the
statement of the Theorem. Let $\mu:(Z,E)\ra (X,D)$ be a
log-resolution which is an isomorphism above $U$, where
$E=\cup_{j\in S'}E_j$ is the inverse image $\mu^{-1}(D)$. By
above, $V^q_i(Z,E)$ is a finite union of sets of the form
$P'\times \cU'$, where $\cup_{P'}P'$ is the canonical
decomposition of $B(Z,E)$ into rational convex polytopes in
$\bR^{S'}$ from Lemma \ref{rebd}, and $\cU'\subset\Pic^\tau(Z,E)$
is a torsion translate of a complex subtorus of $\Pic^\tau(Z)$. By
Lemma \ref{indep}, there is a one-to-one correspondence between
$V^q_i(Z,E)$ and $V^q_i(X,D)$ under the isomorphism of Proposition
\ref{comparison}. The inverse image under this correspondence of
$P'\times\cU'$ is $P\times\cU$ obtained as follows. $P$ is the
image of $P'$ under the projection of $\bR^{S'}$ onto $\bR^S$
given by  the coefficients of the strict transforms of the $D_i$.
Conversely, a realizable boundary $\alpha$ of $X$ on $D$ maps to
the realizable boundary $\{e\}$ of $Z$ on $E$, where $e\cdot
E=\mu^*(\alpha\cdot D)$. Hence $\cup_{P}P$ form a decomposition of
$B(X,D)$ into rational convex polytopes which is a refinement of
the canonical decomposition. Under the isomorphism $\mu^*_{par}$,
subtori of $\Pic^\tau(Z)$ correspond to subtori of $\Pic^\tau(X)$
since $\Pic^\tau(Z)$ consists of finitely many copies of
$\Pic^\tau(X)$. Similarly, torsion elements in $\Pic^\tau(Z,E)$
also correspond to torsion elements in $\Pic^\tau(X,D)$. Therefore
$\cU'$ corresponds to $\cU$ which is a torsion translate of a
subtorus in $\Pic^\tau(X)$. \ \ $\square$

\begin{lem}\label{fin} With the same notation is in Theorem \ref{gv},
fix $q$. Then there exists $i_0$ such that $i>i_0$ implies
$V_i^{q}=\emptyset$.
\end{lem}
\begin{proof} It is enough to restrict to the case when $D$ has
simple normal crossings. For a fixed $k$, it is enough to show
that $V^q_i(L_k)=\emptyset$ for  $i\gg 0$. But $V^q_i(L_k)$ are
closed subsets in the Zariski topology of $\Pic^\tau(X)$ and
$V^q_i(L_k)\supset V^q_{i+1}(L_k)$. Hence $V^q_i(L_k)=\emptyset$
for $i\gg 0$ and this proves the Lemma.\ \
\end{proof}

\begin{example}\label{ex1ct1} (a) Consider the case of Example \ref{ex1}. Let
$\mu:(\wti{\bP},E)\ra (\bP,D)$ be the log resolution given by the
blow-up of the intersection point of the $D_i$'s. Let $E_i$ be the
strict transform of $D_i$ ($i=1,..,5$), and $E_6$ be the
exceptional divisor. The decomposition of $B(\wti{\bP},E)$ is
$\coprod_{k=0,..,4}Q_k$, where $Q_k=[0,1)^6\cap
l_{\wti{\bP}}^{-1}(q_k)$ with
$l_{\wti{\bP}}(e)=(\sum_{i=1,..5}e_i)[\mu^*\cO_{\bP}(1)]+(e_6-\sum_{i=1,..5}e_i)[E_6]$,
$q_k=[L_k]$, $L_k=\mu^*\cO_\bP(k)\otimes\cO_{\wti{\bP}}(-kE_6)$,
for $k=0,..,4$. This induces a decomposition of $B(\bP,D)$ via the
projection $\bR^6\ra \bR^5$ given by $e\mapsto
(e_1,e_2,e_3,e_4,e_5)$. In this case, the induced decomposition on
$B(\bP,D)$ is the same as the canonical one, the image of $Q_k$ is
$P_k$. To check for membership of $P_k$ in $V^q_i$ we need
$h^q(\wti{\bP},\omega_{\wti{\bP}}\otimes L_k)\ge i$. We have
$h^q(\wti{\bP},\omega_{\wti{\bP}}\otimes
L_k)=h^q(\wti{\bP},\mu^*\cO_{\bP}(k-3)\otimes\cO_{\wti{\bP}}((1-k)E_6))$
and we can check that the non-trivial $V^q_i$'s are $V^1_1=P_2\cup
P_3\cup P_4$, $V^1_2=P_3\cup P_4$, $V^1_3=P_4$, $V^2_1=P_0$.

\medskip

(b) This basic example can be computed by other means (see also
Example 1 in Sec. 5 of \cite{Li6}). We included it to show the
basic steps of the method of this article of obtaining the
structure of $V^q_i$: firstly resolve singularities, secondly get
a pool of candidate polytopes from the map $l$, lastly check which
of the line bundles attached to these polytopes satisfy the
conditions of $V^q_i$. This method is different than the method of
\cite{Lib}, \cite{Lib4}, \cite{Li2} which firstly computes the
global polytopes of quasiadjunction and then checks for them the
conditions imposed by $V^q_i$. For example, $P_1$ is not a global
polytope of quasiadjunction since the multiplier ideal sheaf
$\cJ(\alpha\cdot D)=\cO_{\bP}$ for $\alpha\in P_1$. We do not know
how to adapt the method of polytopes of quasiadjunction to prove
Theorem \ref{gv}.

\end{example}

\section{Congruence covers and polynomial periodicity}\label{congcov}

In this section we prove Theorem \ref{pp}. First we reduce the
question  to computing torsion in the canonical stratifications of
Theorem \ref{gv}. Then we use the structure of these sets
described in Theorem \ref{gv} to further reduce the question to
computing lattice points in convex rational polytopes.

\smallskip

\noindent{\bf Reduction to torsion in canonical stratifications.}
With the notation as in Theorem \ref{pp}, denote by $G_N$ the
finite group $H_1(U,\bZ_N)$.
 Denote by $\Pic^\tau(X,D)[N]$ the
$N$-torsion part of $\Pic^\tau(X,D)$ and by $V^q_i$ the sets
$V^q_i(X,D)$. By Corollary \ref{ss3} and Lemma \ref{torsion},
\begin{align*}
h^q(N) & =\sum_{\chi\in G_N^*}h^{n-q}(X,\omega_X\otimes
L_\chi\otimes \cJ(\alpha_\chi\cdot D))\\
& =\sum_{(L,\alpha)\in \Pic^\tau(X,D)[N]}h^{n-q}(X,\omega_X\otimes
L\otimes \cJ(\alpha\cdot D))\\
& = \sum_{i\ge 1}i\cdot \#\left
[(V_i^{n-q}-V_{i+1}^{n-q})[N]\right],
\end{align*}
where $\#S$ denotes the number of elements of a finite set $S$,
and $S[N]$ denotes the set of  $N$-torsion elements of
$\Pic^\tau(X,D)$ lying in the subset $S$. Since
$V_{i+1}^{q}\subset V_i^q$,
$$h^q(N)=\sum_{i\ge 1}\#V_i^{n-q}[N].$$
By Lemma \ref{fin}, this is a finite sum.  By Theorem \ref{gv} and
the inclusion-exclusion formula,
$$h^q(N)=\sum_j a_j\cdot \# [(\cP_j\cdot T_j\cdot \cT _j )[N]],$$
where the sum is finite, $a_j\in \bZ$, $T_j\in \Pic^\tau(X,D)$ is
torsion, $\cT_j$ are subtori of $\Pic^\tau(X)$, and
$$\cP_j=\{ (L_j,\alpha) \in \Pic^\tau(X,D) : \alpha \in P_j
\},$$ for some $L_j$, where $P_j$ is a rational convex polytope
 in $\bR^{S}\cap [0,1)^{S}$. The data $(a_j, L_j, P_j, T_j, \cT_j)$ depends on $q$
but not on $N$. Therefore Theorem \ref{pp} follows from the
assertion that the functions
$$\# [(\cP\cdot T\cdot \cT )[N]]$$
are quasi-polynomials in $N$, where $(\cP, T, \cT)$ is $(\cP_j,
T_j, \cT_j)$ for some $j$. We will show that this assertion is
equivalent to one in terms of lattice points in rational convex
polytopes.

\medskip

\noindent{\bf Reduction to lattice points in polytopes.} Fix $L,
P, T,$ and $\cT$ as above.  The $N$-torsion elements of $\cP\cdot
T\cdot \cT$ consist of realizations of boundaries $(L\otimes
T\otimes M,\alpha)$ in $\Pic^\tau(X,D)$ such that $\alpha\in P$,
$M\in\cT$, $N\alpha\in \bZ^S$, and
$$\cO_X(N\alpha\cdot D)\otimes (L\otimes T)^{\otimes -N}=M^{\otimes
-N}.$$ Since $\cT$ is a torus,
$$\# [(\cP\cdot T\cdot \cT )[N]]=
N^{\text{rank}(\cT)}\cdot\#\cA_N,$$ where
$$\cA_N=\left\{ \alpha\in \bR^S : \alpha \in P,\ N\alpha\in \bZ^S,\
\cO_X(N\alpha\cdot D)\otimes (L\otimes T)^{\otimes -N}\in \cT
\right\}.$$ We will show that $\#\cA_N$ is a quasi-polynomial in
$N$.

\medskip

Let $\beta$ be a point with integral coordinates in the $\bR$-span
of $P$. It follows from Proposition \ref{struc} that we can choose
$L$ such that some multiple of $L$ is linearly equivalent to the
same multiple of $\beta\cdot D$. Let $\cL=\cO_X(L+T-\beta\cdot
D)$, so that $\cL$ is a torsion line bundle. Let $Q$ be the convex
integral polytope $P-\beta$ in $\bR^S$, which might be missing
some of its faces if $P$ does. Via the transformation
$\alpha\mapsto \alpha-\beta$ there is a one-to-one correspondence
of $\cA_N$ with
$$\cB_N=\left \{ \alpha\in \bR^S : \alpha\in Q,\ N\alpha\in\bZ^S,\
\cO_X(N\alpha\cdot D)\in \cL^{\otimes N}\cdot \cT \right \}.$$ We
will show that $\#\cB_N$ is a quasi-polynomial in $N$.

\medskip

Let $W\subset\bR^{S}$ be the subspace consisting of $\alpha$ such
that $\alpha\cdot [D]=0$ in $H^2(X,\bR)$. We can write
$$\cB_N=\frac{1}{N}\rho^{-1}(\cL^{\otimes N}\cdot \cT)\cap W\cap
Q,$$ where $\rho:\bZ^S\ra \Pic(X)$ is the group homomorphism
sending $\alpha$ to $\cO_X(\alpha\cdot D)$. The assertion that
$\#\cB_N$ is a quasi-polynomial follows from Lemma \ref{red}
below. This ends the proof of Theorem \ref{pp}. \ \ $\square$

\medskip

In what follows we do not require that a convex polytope contains
all its faces.

\begin{lem}\label{red} Let $\cG$ be an abelian group and
$\rho:\bZ^S\ra\cG$  a group homomorphism. Let $Q\subset \bR^S$ be
a convex integral polytope. Let $\cT$ be a subgroup of $\cG$. Let
$g\in\cG$ be a torsion point. Let $W$ be a vector subspace of
$\bR^S$. Then, for $N\in \bN$, the function
$$F(N)=\#\left [\frac{1}{N}\rho^{-1}(Ng+\cT)\cap W\cap Q\right ]$$
is a quasi-polynomial.
\end{lem}
\begin{proof}
Denote by $\Lambda$ the inverse image under $\rho$ of $\cT$. Then
$\Lambda$ is a free finitely generated abelian subgroup of
$\bZ^S$. Fix $w\in \rho^{-1}(g+\cT)$. Then for every $N\ge 0$
$$\rho^{-1}(Ng+\cT)=Nw+\Lambda.$$
Moreover, since $g$ is torsion, a multiple of $w$ lies in
$\Lambda$. Let $V=\Lambda\otimes_\bZ\bR$.

\medskip

Let $\Lambda'=\Lambda\cap W$ and $V'=V\cap W$. Thus $\Lambda'$ is
a lattice in $V'$. Fix $w'\in (w+\Lambda)\cap W$. Hence $w'$ is
$\Lambda '$-rational point of $V'$. Then
$$(Nw+\Lambda)\cap W =Nw'+\Lambda ',$$
for $N\ge 0$. Let $Q'=Q\cap V'$, so that $Q$ is a $\Lambda
'$-rational convex polytope in $V'$. Then
\begin{align*}
F(N) &=\#\left[ \frac{1}{N}\left (Nw'+\Lambda'\right)\cap Q'\right
]\\
& =\# \left [\Lambda '\cap N(Q'-w')\right ].
\end{align*}
Since $(Q'-w')$ is a $\Lambda'$-rational convex polytope in $V'$,
 the Lemma follows from Theorem \ref{generh}.\ \
\end{proof}

\begin{thm}\label{generh} (E. Ehrhart, see \cite{St}-4.6.25.)
Let $Q$ be a convex rational polytope in $\bR^n$. Then, for $N\in
\bN$, the function
$$f(N)=\#\left [\bZ^n \cap
NQ\right ]$$ is a quasi-polynomial in $N$.
\end{thm}

The  result above is due to E. Ehrhart and is originally stated
assuming the polytope is closed. The case of polytopes missing
some faces follows since the sum and the difference of two
quasi-polynomials is again a quasi-polynomial.

\section{Generalizations}

We show in this section how the proofs for  Theorem \ref{gv} and
Theorem \ref{pp} extend to give Theorem \ref{gvgen} and Theorem
\ref{ppgen}. Let $U$ be a nonsingular quasi-projective variety.
Let $G$ be a finite abelian group with a surjection $H_1(U)\ra G$.
Let $g:V\ra U$ be the corresponding unramified abelian cover.

\begin{lem}\label{gendecomp} With the notation as above, let $\cW$
be a unitary local system on $U$.

\noindent (a) There is an eigenspace decomposition
$$H^m(V,g^*\cW)=\bigoplus_{\chi\in G^*} H^m(U,\cW\otimes \cV_\chi),$$
where $\cV_\chi$ is the rank one unitary local system induced by
$\chi$.

\noindent (b) The decomposition above is compatible with the Hodge
filtration, i.e. there is an eigenspace decomposition

$$\Gr_F^pH^m(V,g^*\cW)=\bigoplus_{\chi\in
G^*}\Gr_F^pH^m(U,\cV_\chi\otimes\cW).$$
\end{lem}
\begin{proof} (a) Since $g$ is finite,
$H^m(V,g^*\cW)=H^m(U,g_*g^*\cW)$ which, by the projection formula,
is isomorphic to $H^m(U,\cW\otimes g_*\bC_V)$. The claim follows
since $g_*\bC_V=\oplus_{\chi\in G^*}\cV_\chi$ (see Remark
\ref{ggg}).

\medskip

\noindent (b) Follows from (a) by the $G$-equivariant version of
functoriality for the mixed Hodge structure on the cohomology of
unitary local systems (see \cite{Ti}).
\end{proof}

Let $X$ be a nonsingular compactification of $U$ with complement
$D=\cup_{i\in S}D_i$ a divisor with simple normal crossings. For a
unitary local system $\cW$ on $U$, denote from now on by
$\overline{\cW}$ the vector bundle of the canonical Deligne
extension of $\cW$ to $X$.

\begin{rem}\label{canonical extension rank one} (a) Let $\cV$ be a unitary local system of rank one on $U$ and let
$(L,\alpha)\in\Pic^\tau(X,D)$ be the corresponding realization of
boundary. The inverse of $(L,\alpha)$ in $\Pic^\tau(X,D)$ is
$(M,\beta)$, where $M=L^{-1}+\sum_{\alpha _i\ne 0}D_i$, and $\beta
_i$ is $0$ if $\alpha _i=0$ and $1-\alpha _i$ otherwise. By the
proof of Theorem \ref{ls}, $\overline{\cV}=M^{-1}$, hence
$\overline{\cV}=L-\sum_{\alpha _i\ne 0}D_i.$

\medskip

(b) In general, the canonical extension is not compatible with $\otimes$. For example, if $\cV$ and $\cW$ are two rank one unitary local systems on $U$ with corresponding realization of boundaries $(L,\alpha)$, respectively $(M,\beta)$, then
\begin{align*}
\overline{\cV\otimes\cW} &=L+M-\rndown{\alpha +\beta}\cdot D-\sum_{\{\alpha _i+\beta _i\}\ne 0}D_i,\\
\overline{\cV}\otimes\overline{\cW} &=L+M-\sum_{\alpha _i\ne 0}D_i-\sum_{\beta _i\ne 0}D_i.
\end{align*}
However, if $\cW$ is the restriction to $U$ of a unitary local
system (of arbitrary rank) on $X$, then
$\overline{\cV\otimes\cW}=\overline{\cV}\otimes\overline{\cW}$.
This follows from the construction of the canonical extension
\cite{De}.

\medskip

(c) If $\cV$ and $\cW$ are unitary local systems on $U$, with
$\cV$ of rank one corresponding to $(L,\alpha)$ and $\cW$ being
the restriction of a local system from $X$, by (a) and (b) above
$$\overline{\cV\otimes\cW}=L\otimes\cO_X(-\sum_{\alpha _i\ne
0}D_i)\otimes\overline{\cW}.$$ The last vector bundle determines
polytopes in $[0,1)^S$ such that the vector bundle remains
constant when $\alpha$ varies within these polytopes. For example,
consider the polytopes $P_k$ of (\ref{desc}). Then
$\overline{\cV\otimes\cW}$ remains constant if $\alpha$ varies
within a fixed polytope $P_{k,S'}:=P_k\cap \{ \alpha\ |\ \alpha
_i\ne 0 \text { if }i\in S', \alpha _i=0 \text{ if }i\notin S'\}$,
where $S'\subset S$. The polytopes $P_{k,S'}$ form a more refined
decomposition of $B(X,D)$ than (\ref{desc}).

\medskip

(d) A similar conclusion holds for any $\cV$ and $\cW$ as in (c)
without assuming that $\cW$ is the restriction of a local system
from $X$. More precisely, there exists a finer decomposition than
(\ref{desc}) of $B(X,D)$ into rational convex polytopes $P_k$ such
that if $\cV=(L,\alpha)$ varies with $\alpha\in P_k$ and $L$ and
$k$ fixed, then $\overline{\cV\otimes\cW}$ is constant. This
follows from the explicit construction of the canonical
extensions. For example, let $u_i$ be a frame of the multi-valued
flat sections of $\cV$ around a general point of $D_i$. Then
$\overline{\cV}$ is locally given by the holomorphic frame
$v_i=\exp(\log z_i\cdot\alpha _i)u_i$, where $z_i$ is a local
equation for $D_i$ and the monodromy of $\cV$ around $D_i$ is
$\exp(-2\pi i\alpha _i)$ with $\alpha_i\in [0,1)$. Fix an
orthonormal frame $u_{i,j}$ $(1\le j\le \rank (\cW))$ of the
multi-valued flat sections of $\cW$ around a general point of
$D_i$, such that the monodromy around $D_i$ sends $u_{i,j}$ to
$\exp(-2\pi i\alpha_{i,j})u_{i,j}$, with $\alpha_{i,j}\in [0,1)$.
Then a local holomorphic frame of $\overline{\cW}$ is given by
$v_{i,j}:=\exp(\log z_i\cdot \alpha _{i,j}) u_{i,j}$. A local
holomorphic frame of $\overline{\cV\otimes\cW}$ is given by
$w_{i,j}:=\exp(\log z_i\cdot \beta _{i,j}) u_i\otimes u_{i,j}$,
where $\beta _{i,j}=\{\alpha_i+\alpha_{i,j}\}$ is the fractional
part of $\alpha_i+\alpha_{i,j}$. The polytopes are determined by
imposing the condition that $\beta$ is constant (here
$\alpha_{i,j}$ are fixed). Remark that
$\overline{\cV}\otimes\overline{\cW}$ is given locally by
$v_i\otimes v_{i,j}=\exp(\log
z_i\cdot\rndown{\alpha_i+\alpha_{i,j}})w_{i,j}$.

\end{rem}

\noindent {\bf Proof of Theorem \ref{gvgen}.} The proof is
essentially the same as for Theorem \ref{gv}. Assume first that
$D$ is a divisor with simple normal crossings. Then, by the
decomposition (\ref{desc}),
$$W^{p,q}_i(U,\cW)=\coprod _k W^{p,q}_i(U,\cW)\cap [P_k\times \Pic^\tau(X)].$$
The intersection is possible via the identification of unitary local systems of rank one on $U$ with realizations of boundaries of $X$ on $D$. Pointwise, the set $W^{p,q}_i(U,\cW)\cap [P_k\times \Pic^\tau(X)]$ consists of the local systems $\cV\in\Hom(H_1(U),S^1)$ corresponding to $\{(L_k+M,\alpha)\ |\ \alpha\in P_k, M\in\Pic^\tau(X)\}$ such that $\dim \Gr^p_FH^{p+q}(U,\cV\otimes\cW)\ge i$, where $L_k$ depends on $P_k$. By \cite{Ti}-2nd Theorem, part (a),
$$\dim \Gr_F^pH^{p+q}(U,\cV\otimes\cW)=h^q(X,\Omega_X^p(\log D)\otimes \overline {\cV\otimes\cW}).$$
By Remark \ref{canonical extension rank one} -(d), the last
quantity is constant if $\alpha$ varies within a fixed polytope
$P_k$, possibly coming from a finer decomposition of $B(X,D)$ than
(\ref{desc}). We will work with this finer decomposition from now
on and use the same notation, $P_k$, for its polytopes.

\medskip

 Fix $\alpha_{k}\in
P_{k}$ and let $\cV_{k}$ denote the local system corresponding to
$(L_k,\alpha_{k})$. The proof of Proposition \ref{struc} goes
word-by-word for the finer $P_k$'s. In other words, we can assume
that $(L_k,\alpha_{k})$ is a torsion element of $\Pic^\tau(X,D)$.
We have
$$W^{p,q}_i(U,\cW)=\coprod_{k} P_{k}\times W^{p,q}_i(\cV_{k}\otimes\cW),$$
where
$$W^{p,q}_i(\cV_{k}\otimes\cW):=\{\
\cV\in\Hom(H_1(X),S^1)\ |\
\dim\Gr_F^pH^{p+q}(U,\cV_{k}\otimes\cW\otimes\cV_{|U})\ge i\ \}.$$
Since $\cV_{k}$ have finite order, we can find a finite abelian
group $G$ with and embedding
$G^*\subset\Pic^\tau(X,D)\cong\Hom(H_1(U),S^1)$, such that
$\cV_{k}\in G^*$. Let $\chi_{k}$ be the character of $G$ given by
$\cV_{k}$. Let $g:V\ra U$ the unramified $G$-cover given by
$G^*\subset\Pic^\tau(X,D)$. Then, by Lemma \ref{gendecomp}-(b),
$$W^{p,q}_i(\cV_{k}\otimes\cW)=\{ \cV\in\Hom(H_1(X),S^1)\ |\
\dim\Gr_F^pH^{p+q}(V,g^*(\cW\otimes\cV_{|U}))_{\chi_{k}}\ge i\
\}.$$

If $H_1(X)=0$, then $\Pic^\tau(X)=\{(\cO_X,0)\}$, $\Pic^\tau(X,D)=B(X,D)$ and $W^{p,q}_i(\cV_{k,S'}\otimes\cW)$ is either empty or $\{(\cO_X,0)\}$. Then the statement of Theorem \ref{gvgen} for this case is clear.

\medskip

If $\cW$ is the restriction of a local system from $X$,
then the
statement of Theorem \ref{gvgen} for this case follows from
Theorem \ref{sigen}.

\medskip

If $D$ is assumed to be general, the conclusion follows from the
simple normal crossings case as in the proof of Theorem \ref{gv}.\
\ $\square$

\medskip

Let $L(X)$ denote the category of local systems on $X$, $M(X)$
denote the Betti realization of the moduli space of local systems
on $X$, and $M_i(X)$ denote the component of local systems of rank
$i$. This changes the notation from the Introduction. Thus
$M_1(X)$ is now what we denoted by $M_B(X)$ in the Introduction.
Let $U_1(X)$ denote the subspace of unitary local systems of
$M_1(X)$. Note that in the notation of the following result, $\cW$
will be a local system on $X$, not on $U$.

\begin{thm}\label{sigen}  With the notation as in Theorem
\ref{si}, let $U$ be an open subset of $X$ and such that $D=X-U$
is a simple normal crossings divisor. Let $V=f^{-1}(U)$. Let $\cW$
be a unitary local system on $X$.
Then for all $\chi\in G^*$, the
image of the set
\begin{equation}\label{loci_with_Gr_and_chi}
\{ \cV\in U_1(X) \ |\
\dim\Gr_F^pH^{p+q}(V,f^*(\cW\otimes\cV)_{|V})_\chi \ge i \}
\end{equation}
in $\Pic^\tau(X)$ is a finite union of torsion translates of
complex subtori of $\Pic^\tau(X)$, and so is any intersection of
these translates.
\end{thm}

\begin{proof}
As for Theorem \ref{si}, this  follows from the general theory of
absolute constructible sets of C. Simpson \cite{Si}. Any absolute
closed subset of $M_1(X)$ is a finite union of torsion translates
of triple subtori of $M_1(X)$ (\cite{Si}-Theorem 6.1). A triple
subtorus is a closed connected real analytic subgroup of the
underlying analytic group of $M_1(X)$ such that it is an algebraic
subgroup in each three complex structures Betti, de Rham, and
Dolbeault (see Introduction). An intersection of absolute
constructible subsets is also absolute. A way of producing
absolute sets is via absolute functors between two categories of
local systems of smooth projective varieties and via absolute
natural transformations between absolute functors. For example,
the image and the inverse image of absolute sets under absolute
functors is again absolute (\cite{Si}-Corollary 7.3 and Lemma
7.4).

\medskip

First, consider the loci
\begin{equation}\label{loci_in_M(X)}
\{\ \cV\in M(X)\ |\ \dim H^m(U,\cV_{|U})\ge i\ \},
\end{equation}
where $m$ and $i$ are fixed. In case $U=X$, these loci are
absolute closed sets (\cite{Si}-Corollary 7.14), due to the fact
that $H^m(X,.)$ defines an absolute functor from $L(X)$ to
$L(\text{point})$. We will reduce the general case to this
particular case. Recall that the loci (\ref{loci_in_M(X)}) are
closed subsets being inverse images under the restriction morphism
$M(X)\ra M(U)$ of closed subsets (see \cite{Ar2}- Corollary 2.5).

\medskip

We will not be able to prove that $H^m(U,\cV_{|U})$ is an absolute
functor from $L(X)$ to $L(\text{point})$. However, for the loci
given by the dimension of $H^m(U,\cV_{|U})$ to be absolute, it is
enough to prove that there exists a filtration of
$H^m(U,\cV_{|U})$ such that the graded pieces form absolute
functors. Indeed, the associated graded functor will be an
absolute functor to vector spaces of same dimensions as given by
$H^m(U,\cV_{|U})$.

\medskip

First, by \cite{Di}-2.4, there exists a long exact sequence
\begin{equation}\label{exact_spectral}
\ldots\ra H^m_D(X,\cV)\xrightarrow{\phi^{m}=\phi^m(\cV)}
H^m(X,\cV)\ra H^m(U,\cV)\ra H^{m+1}_D(X,\cV)\ra\ldots.
\end{equation} Let $i:D\hookrightarrow X$ denote the inclusion
map. Here $H^m_D(X,\cV)=\bH^m(D,i^!\cV)$. By \cite{Di}-Proposition
3.3.7, $i^!\cV=\bD_D(\cV^\vee_{|D})[-2n]$, where $\bD_D(.)$ is the
dual complex in the bounded derived category of constructible
sheaves on $D$, and $n$ is the (complex) dimension of $X$. Thus,
$H^m_D(X,\cV)=\bH^{m-2n}(D,\bD(\cV ^\vee_{|D}))$, and by
Poincar\'{e}-Verdier duality (\cite{Di}-Theorem 3.3.10), this
space is naturally isomorphic with
$H^{2n-m}(D,\cV^\vee_{|D})^\vee$. We will use a Mayer-Vietoris
spectral sequence degenerating to this last space, whose
$E_\infty$ term is an absolute functor and induces filtrations on
$\ker \phi^m$ and $\text{coker}\, \phi^m$ such that the associated
graded functors are absolute functors. This will suffice for our
purposes since
$$\dim H^m(U,\cV_{|U})=\dim [\text{coker}\,\phi^{m}\oplus\ker
\phi^{m+1}].$$

\medskip

Let $\psi^p$ be the restriction map $H^p(X,\cV)\ra H^p(D,\cV)$.
The natural transformations $\phi$ and $\psi$ are related by
$\phi^m(\cV)^\vee=\psi^{2n-m}(\cV^\vee)$. We discuss first $\psi$
and then we draw conclusions about $\phi$.

\medskip

Consider the Mayer-Vietoris spectral sequence
\begin{equation}\label{MVspectralseq}
E_1^{p,q}(D)=H^q(D^{[p]},\cV)\Rightarrow H^{p+q}(D,\cV),
\end{equation}
where $D^{[p]}$ is the disjoint union of $p+1$-fold intersections
of distinct irreducible components of $D=\cup_{i\in S}D_i$. For
$I\subset S$, let $D_I=\cap_{i\in I}D_i$. Then
$D^{[p]}=\coprod_{I\subset S; \#I=p+1}D_I$. We recall next the
construction of (\ref{MVspectralseq}). Define
$$C^{p,q}(D):=\mathop{\bigoplus_{I\subset
S}}_{\#I=p+1}C^q(D_I,\cV),$$ where $C^q(D_I,\cV)$ denote the
$\cV$-twisted cochains. Let $\hat{d}:C^{p,q}(D)\ra C^{p+1,q}$ and
$d: C^{p,q}(D)\ra C^{p,q+1}(D)$ denote the  Cech and cochain
differentials. Then $(C^{\bullet,\bullet}(D),\hat{d},d)$ becomes a
double complex. We suppressed from notation that
$C^{\bullet,\bullet}(D)$ depends on the closed cover
$\{D_i\}_{i\in S}$ and the local system $\cV$. Let
$Tot^m(D)=\oplus_{p+q=m}C^{p,q}(D)$, so that
$(Tot^\bullet(D),\hat{d}+d)$ is the total complex of
$(C^{\bullet,\bullet}(D),\hat{d},d)$. Then (\ref{MVspectralseq})
is the spectral sequence $(E^{\bullet,\bullet}_r(D),d_r)$ given by
the filtration of $Tot^\bullet(D)$ defined by $$F^p
Tot^m(D)=\mathop{\bigoplus_{r+s=m}}_{r\ge p}C^{r,s}.$$

\medskip

Define $C^{p,q}(X)$ to be $0$ if $p\ne 0$ and $C^q(X,\cV)$ if
$p=0$. Let $\hat{d}:C^{p,q}(X)\ra C^{p+1,q}$ and $d: C^{p,q}(X)\ra
C^{p,q+1}(X)$ denote the  $0$-map and the cochain differential.
Then $(C^{\bullet,\bullet}(D),\hat{d},d)$ becomes a double
complex. This construction is similar to $C^{p,q}(D)$ in the sense
that it depends on the closed cover $\{X\}$ of $X$ and the local
system $\cV$; the Cech differential is $0$ for this closed cover.
The corresponding spectral sequence, denoted by
$(E^{\bullet,\bullet}_r(X),d_r)$, degenerates to $H^{p+q}(X,\cV)$.
Also, $E_r^{p,q}(X)=H^q(X,\cV)$ if $p=0$ and equals $0$ otherwise,
for $r\ge 1$.

\medskip

Define a map of double complexes of bidegree $(0,0)$
$$\psi^{\bullet,\bullet}:(C^{\bullet,\bullet}(X),\hat{d},d)\ra (C^{\bullet,\bullet}(D),\hat{d},d)$$
as follows. For $p\ne 0$, $\psi^{p,q}=0$. Define
$\psi^{0,q}:C^q(X,\cV)\ra \oplus_{i\in S} C^{q}(D_i,\cV)$ as the
direct sum of the restriction maps. It is trivial that
$\psi^{\bullet,\bullet}$ commutes with $d$ ($0$ goes to $0$). Note
that $\psi^{0,q}$ is the composition of the restriction map
$C^q(X,\cV)\ra C^q(D,\cV)$ with the Cech map $C^q(D,\cV)\ra
\oplus_{i\in S} C^{q}(D_i,\cV)$. Hence $\psi^{\bullet,\bullet}$
commutes with $\hat{d}$ also.

\medskip

The map of double complexes $\psi^{\bullet,\bullet}$ induces of
morphism of filtered total complexes (which we prefer not to have
a notation for, so that we do not confuse it with the cohomology
restriction map $\psi^p$). By \cite{Mc} -Theorem 3.5, we have a
morphism of spectral sequences
$$\psi_r^{\bullet,\bullet}: (E^{\bullet,\bullet}_r (X),d_r)\ra
(E^{\bullet,\bullet}_r(D),d_r),\ \ \ \ \  (0\le r\le \infty).$$ In
particular, we have a map
$$
\psi_\infty ^m : E_\infty ^m (X)= H^{m}(X,\cV) \ra E_\infty
^{m}(D)=\bigoplus _{p+q=m} E_\infty ^{p,q}(D),
$$
where $E_\infty ^{p,q}(D)=\Gr_L^p H^{p+q}(D,\cV)$ for some
filtration $L$, and $E_\infty ^{m}(D)$ is the associated graded
vector space of $H^m(D,\cV)$.

(Example: If $D$ has only two components $D_1$ and $D_2$, then
$E_\infty ^{m}(D)=\ker {{d}}\; _1 ^{0,m} \oplus \text{coker}\;
{{d}}\; _1^{0,m-1}$, where $d\; _1^{0,m}:H^m(D_1,\cV)\oplus
H^m(D_2,\cV)\ra H^{2}(D_1\cap D_2,\cV)$ is the map induced by the
Cech differential $\hat {d}$. Then $\psi_\infty ^m:H^m(X,\cV)\ra
E_\infty ^m(D)$ is the map $(\psi_1^{0,m},0)$, where
$\psi_1^{0,m}:H^m(X,\cV)\ra H^m(D_1,\cV)\oplus H^m(D_2,\cV)$ is
the direct sum of the restriction maps.)

\medskip

Since $\psi^{\bullet,\bullet}$ factors through the restriction map
$C^q(X,\cV)\ra C^q(D,\cV)$, the map $\psi_\infty ^m$ factors as
$\Gr_L \circ \psi^m$, where $\psi^m: H^m(X,\cV)\ra H^m(D,\cV)$ is
the cohomology restriction map and $\Gr_L: H^m(D,\cV)\ra E_\infty
^{m}(D)$ is taking elements to their associated graded image. Thus
$\ker$ and $\text{coker}$ of $\psi_\infty^m$ have the same
dimension as $\ker$ and, respectively, $\text{coker}$ of $\psi^p$.

\medskip

The $D_I$ are smooth projective and the differential $d_1$ of
(\ref{MVspectralseq}) is obtained from restrictions. It follows
that the $E_1(D)$ terms of (\ref{MVspectralseq}) form absolute
functors from local systems on $X$ to the point, and the
differential $d_1$ is an absolute natural transformation (see
\cite{Si}-Proposition 7.8 and 7.9). Kernels and cokernels of
absolute natural transformations give absolute functors
(\cite{Si}-Lemma 7.13), hence $(E_r(D),d_r)$ is an absolute
functor with an absolute natural transformation, and so
$E_\infty(D)$ is an absolute functor. Also $\psi_\infty ^p$ is an
absolute natural transformation since it is obtained from
restriction maps. Thus $\ker \psi_\infty^p$ and
$\text{coker}\,\psi_\infty^p$ are absolute functors.

\medskip

By dualizing the terms in the spectral sequences we obtain a
natural transformations $\phi_\infty^m$ such that
$\phi_\infty^m(\cV)^\vee=\psi_\infty^{2n-m}(\cV^\vee)$. We want to
show that $\phi_\infty^m$ are also absolute. It is enough to prove
that  $(E_1(D),d_1)$ dualizes to an absolute functor with an
absolute natural transformation. Since composition of absolute
natural transformation is again absolute, it is enough then to
show that the duality $H^m(X,\cV)\longleftrightarrow
H^{2n-m}(X,\cV^\vee)^\vee$ (in our case we will want to replace
$X$ by any $D_I$) is an absolute natural transformation in both
directions. The dual $\cV\longleftrightarrow\cV^\vee$ is an
absolute functor (\cite{Si}-Lemma 7.11), so $\cV\mapsto
(\cV,\cV^\vee)$ is an absolute functor from $L(X)$ to $L(X)\times
L(X)$. Thus the natural transformation $H^m(X,\cV)\otimes
H^{2n-m}(X,\cV^\vee)\ra H^{2n}(X,\bC_X)$ factors through the
absolute natural transformation of cup
product(\cite{Si}-Proposition 7.12)
$$H^m(X,\cV_1)\otimes H^{2n-m}(X,\cV_2)\ra
H^{2n}(X,\cV_1\otimes\cV_2)$$ via an absolute functor, hence it is
itself absolute. Therefore $\phi_\infty^m$, and hence $\ker
\phi_\infty^m$ and $\text{coker}\,\phi_\infty^m$, are absolute
natural transformations.

\medskip

Thus the loci given by the dimension of
$\text{coker}\,\phi_\infty^{m}(\cV)\oplus\ker
\phi_\infty^{m+1}(\cV)$ are absolute sets. This finishes the proof
of the fact that the loci (\ref{loci_in_M(X)}) are absolute closed
sets. By \cite{Si} -Corollary 6.2 and Lemma 7.1 (see also
Conclusion in loc. cit.), the restriction of the loci
(\ref{loci_in_M(X)}) to $M_1(X)$ are finite unions of torsion
translates of subtori of $M_1(X)$.

\medskip

Now, let $\cW$ be a local system on $X$. The same proof as above
extends to show that the loci
\begin{equation}\label{loci_with_W}
\{ \cV\in M(X)\ |\ \dim H^m(U,(\cW\otimes\cV)_{|U})\ge i\ \}
\end{equation}
($m$, $i$ are fixed) are absolute closed subsets. Hence, as above,
the restriction of the loci (\ref{loci_with_W}) to $M_1(X)$ are
finite unions of torsion translates of subtori of $M_1(X)$.

\medskip

With $f$ and $V$ as in the Theorem, the loci
\begin{equation}\label{loci_with_f}
\{ \cV\in M(X)\ |\ \dim H^m(V,f^*(\cW\otimes\cV)_{|V})\ge i\ \}
\end{equation}
are absolute closed sets since they are the intersection of the
loci of local systems $\cV'$ on $Y$ given by
$H^m(V,(f^*\cW\otimes\cV')_{|V})$, which are absolute by the
previous discussion, with the image of the absolute functor $f^*$
(\cite{Si}-Proposition 7.8). Hence, the restriction of the loci
(\ref{loci_with_f}) to $M_1(X)$ are finite unions of torsion
translates of subtori of $M_1(X)$.

\medskip

Assuming now for the first time that $\cW$ is unitary, we prove
the statement of the Theorem for the loci
\begin{equation}\label{loci_with_Gr}
\{\ \cV\in U_1(X)\ |\
\dim\Gr_F^pH^{m}(V,f^*(\cW\otimes\cV)_{|V})\ge i\ \}.
\end{equation}
By \cite{Ti}, for $\cV$ a unitary local system on $X$, the
dimension of $\Gr_F^pH^{m}(V,f^*(\cW\otimes\cV)_{|V})$ is computed
from the $m-p$ cohomology of the $p$-th order logarithmic
differentials twisted by the canonical extension of
$f^*(\cW\otimes\cV)_{|V}$ to some nonsingular compactification of
$V$ with boundary a simple normal crossings divisor. Thus the
dimension of $\Gr_F^pH^m(V,f^*(\cW\otimes\cV)_{|V})$ is a
semicontinuous function in $\cV$. Hence, as in \cite{Si}-section 5
(or \cite{Ar1}, see also the proof of Theorem \ref{si}), the
restrictions of the loci (\ref{loci_with_f}) to $U_1(X)$ give a
stratification which is a refinement of the stratification given
by the loci (\ref{loci_with_Gr}). Subtori of $M_1(X)$ restrict to
subtori of $U_1(X)$. Thus the loci (\ref{loci_with_Gr}) are finite
unions of torsion translates of subtori of $U_1(X)$.

\medskip

Lastly, the statement of the Theorem is the $G$-equivariant
version of the above arguments (see Theorem \ref{si}).
\end{proof}

\begin{rem}\label{on_IH} If one only assumes that $\cW$ is a unitary local
system on $U$ in Theorem \ref{sigen}, then the analogue of
(\ref{exact_spectral}) would involve the intersection cohomology
$IH^m(X,\cW\otimes\cV_{|U})$. However we do not know how to prove
that this forms an absolute functor from $L(X)$ to
$L(\text{point})$. One is then led naturally to Question
\ref{qgen} of the Introduction.
\end{rem}

\noindent{\bf Proof of Corollary \ref{corollary_to_general case}.} The loci of the statement can be rewritten as
\begin{equation}\label{union}
\bigcup _{
\begin{array}{c}
(i_0,\ldots, i_m)\\
 i_0+\ldots +i_m\ge i
\end{array}
}\bigcap _{0\le p\le m}W^{p,m-p}_{i_p}(U,\cW).\end{equation}
The analog of Lemma \ref{fin} holds for $W^{p,q}_i(U,\cW)$ too, using the sets $W^{p,q}_i(\cV_k\otimes\cW)$ instead of the sets $V^q_i(L_k)$. That is, the union in (\ref{union}) is finite. Hence the Theorem follows from Theorem \ref{gvgen}.\ \ $\square$

\noindent {\bf Proof of Theorem \ref{ppgen}.} With the notation as
in the proof of Theorem \ref{pp},
$$h^{p,q}_\cW (N)=\sum_{i\ge 1}i\cdot \#\left
[(W_i^{p,q}(U,\cW)-W_{i+1}^{p,q}(U,\cW))[N]\right],$$ by Lemma
\ref{gendecomp} and Lemma \ref{torsion}. The rest of the proof is
word by word as that of Theorem \ref{pp}, with the following
exceptions: $V^q_i$ replaced by $W^{p,q}_i(U,\cW)$; Lemma
\ref{fin} holds for this case as well;  the reference to Theorem
\ref{gv} is replaced by the Theorem \ref{gvgen} for
$W^{p,q}_i(U,\cW)$; and Proposition \ref{struc} works as well for
any refinement of the original decomposition of $B(X,D)$. \ \
$\square$

\end{document}